\documentclass[12pt]{amsart}

\usepackage{amscd}
\newtheorem{theorem}{Theorem}[section]
\newtheorem{lemma}[theorem]{Lemma}

\theoremstyle{definition}

\theoremstyle{remark}
\newtheorem{remark}[theorem]{Remark}
\numberwithin{equation}{section}

\def\A{\bar A}
\def\B{\bar B}
\def\E{\mathcal E}
\def\CB{\bar C}

\def\R{\bar R_0(\lambda)}
\def\fh{\mathfrak h}
\def\fht{\mathfrak h_\Theta}
\def\H{\mathcal H}
\def\HB{\bar{\mathcal H}}
\def\HZ{\mathcal H_0}
\def\HU{\mathcal H_1}
\def\HD{\mathcal H_2}
\def\KU{\mathcal K_1}

\def\HUB{\bar{\mathcal H}_1}
\def\HDB{\bar{\mathcal H}_2}
\def\KUB{\bar{\mathcal K}_1}
\def\KDB{\bar{\mathcal K}_2}
\def\BO{\mathsf B}

\def\vp{\varphi}
\def\uno{\mathsf 1}
\def\taub{\bar\tau}
\def\Gb{{\breve G}(\lambda)}
\def\GCb{{\breve G_C}(\lambda)}
\def\G{G(\lambda)}
\def\GC{{G_C}(\lambda)}
\def\GB{{\bar G}(\lambda)}

\def\GBb{\breve{\bar G}(\lambda)}
\def\taut{\tau_\Theta}
\def\GTb{\breve G_\Theta(\lambda)}
\def\GT{G_\Theta(\lambda)}
\def\GMT{\Gamma_\Theta(\lambda)}

\def\ran{\text{\rm Ran}}
\def\ker{\text{\rm Ker}}
\def\CUB{\bar C_1}
\def\CDB{\bar C_2}
\def\RB{\bar R(\lambda)}
\def\z*{\bar z}

\def\Q{{\mathcal Q}_\Theta}

\def\K{\mathcal K}
\def\V{\mathcal V}

\def\N{\mathcal N}
\def\g{{\mathcal G}}
\def\RE{\mathbb R}
\def\C{{\mathbb C}}

\def\p{\par\noindent}

\def\As{\mathsf A}
\def\Es{\mathsf E}
\def\Phs{\mathsf \Phi}
\def\Pss{\mathsf \Psi}
\def\Pis{\mathsf \Pi}

\begin{document}

\title[Singular Perturbations of Abstract Wave Equations]
{Singular Perturbations of Abstract Wave
Equations
}

\author{Andrea Posilicano}

\address{Dipartimento di Scienze Fisiche e Matematiche, 
Universit\`a dell'Insubria, I-22100
Como, Italy}

\email{posilicano@uninsubria.it}

\begin{abstract}
Given, on the Hilbert space $\H_0$, the self-adjoint
operator $B$ and the skew-adjoint operators $C_1$ and 
$C_2$, we consider, on the Hilbert space $\H\simeq D(B)\oplus\H_0$, 
the skew-adjoint operator $$W=\left[\begin{matrix}
C_2&\uno\\-B^2&C_1\end{matrix}\right]$$ corresponding to the
abstract wave equation
$\ddot\phi-(C_1+C_2)\dot\phi=-(B^2+C_1C_2)\phi$. 
Given then an auxiliary Hilbert
space $\fh$ and a linear map 
$\tau:D(B^2)\to\fh$ with a kernel $\K$ dense in $\H_0$, 
we explicitly construct
skew-adjoint operators $W_\Theta$
on a Hilbert space $\H_\Theta\simeq D(B)\oplus\H_0\oplus \fh$ which coincide with $W$ on
$\N\simeq\K\oplus D(B)$. 
The extension parameter $\Theta$ ranges over the set of positive, bounded and
injective self-adjoint operators on $\fh$. \par
In the case $C_1=C_2=0$ our construction allows a natural definition
of negative (strongly) singular perturbations 
$A_\Theta$ of $A:=-B^2$ such that the diagram 
$$
\begin{CD}
W @>>> W_\Theta\\
@AAA   @VVV\\
A@>>> A_\Theta
\end{CD}
$$
is commutative. 
\end{abstract}

\maketitle

\section{Introduction}
Given a negative and injective self-adjoint operator $A=-B^2$ on the 
Hilbert space $\H_0$ with scalar product
$\langle\cdot,\cdot\rangle_0$ and corresponding norm $\|\cdot\|_0\,$,
we consider the abstract wave equation
$$\ddot\phi=A\phi\,.$$ The Cauchy problem for such an equation is
well-posed and 
$$\phi(t):=\cos tB\,\phi_0+B^{-1}\sin tB\,\dot\phi_0$$ is the (weak)
solution with initial data $\phi_0\in D(B)$ and $\dot\phi_0\in\H$. 
More precisely, using a block matrix operator notation, 
$$
\left[\begin{matrix}
\cos tB&B^{-1}\sin tB\\-B\sin tB&\cos tB \end{matrix}\right]
$$
defines a strongly continuous group of evolution on the Hilbert space 
$\H_1\oplus\H_0$, where $\H_1$
denotes $D(B)$ endowed with the scalar product giving rise to the graph
norm. It preserves
the energy $${\mathcal E(\phi,\dot\phi)}
:=\frac{1}{2}\,\left(\|\dot\phi\|_0^2+\|B\phi\|^2_0\right)$$ and, 
in the case the
Hilbert space $\H_0$ is real, constitutes a
group of canonical transformations with respect to the standard
symplectic form 
$$\Omega((\phi_1,\dot\phi_1),(\phi_2,\dot\phi_2))
:=\langle\phi_1,\dot\phi_2\rangle_0-\langle\phi_2,\dot\phi_1\rangle_0\,.$$
Its generator is given
by $$\overset\circ W=\left[\begin{matrix}
0&\uno\\-B^2&0\end{matrix}\right]\,:\,D(B^2)\oplus D(B)\subseteq
\H_1\oplus\H_0\to \H_1\oplus\H_0\,;$$
it is the Hamiltonian vector field
corresponding, via $\Omega$,
to the Hamiltonian function $\mathcal E$. \par
From the point of view of Hamiltonian systems (with infinite degrees
of freedom) a
more suitable phase space is given by the space of
finite energy states, i.e. the maximal domain of
definition of the
energy $\mathcal E$. This set is given by $D(\E)=\HUB\oplus\HZ$ where $\HUB$
denotes the Hilbert space obtained by completing $D(B)$ endowed with
the scalar product
$$
[\phi_1,\phi_2]_1:=
\langle B\phi_1,B\phi_2\rangle_0\,.
$$
By our injectivity hypothesis $0\notin\sigma_{pp}(A)$, but
$0\in\sigma(A)\backslash\sigma_{pp}(A)$ is not excluded (e.g. when 
$A=\Delta$ and  $\H_0=L^2(\RE^d)$). Thus in general $\HUB$ is not 
contained into $\HZ$.\par
It is then possible to define a new operator $W$ which is proven to 
be skew-adjoint on the Hilbert space $D(\E)$. Such an
operator is nothing but the closure of $\overset\circ W$, 
now viewed as an
operator on the larger space $D(\E)$. By Stone's theorem $W$
generates a strongly continuous group $U^t$ of unitary operators 
which preserves the energy, which now coincides with the norm of 
the ambient space. 
\par 
Consider now a self-adjoint operator $\hat A\not=A$ which is a 
singular perturbation of $A$, i.e the
set $\K:=\{\phi\in D(A)\cap D(\hat A)\,:\, A\phi=\hat A\phi\}$ is dense
in $\HZ$ (see e.g. \cite{[Ko]}). Since $\K$ is closed with respect to the graph norm on
$D(A)$, the linear operator $A_\K$, obtained by restricting $A$ to the
set $\K$, is a densely defined closed symmetric operator. Therefore
the study of singular perturbations of $A$ is brought back to the
study of self-adjoint extensions of the symmetric operators obtained
by restricting $A$ to some dense, closed with respect to the graph
norm, set. We refer to \cite{[AK]} and its huge list of references for
the vast literature on the subject. However here we found more convenient
to use the approach introduced in \cite{[P1]}.\par 
In the case the singular perturbation $\hat A$ is negative and injective, 
we are interested in
describing $\hat W$, the analog
of $W$ relative to $\hat A$. A natural question is:\p
1. Is $\hat W$ a singular perturbation of $W$?\p
Here a skew-adjoint 
operator $\hat W$ on $D(\hat\E)\supseteq D(\E)$ is said
to be a singular perturbation of the skew-adjoint operator $W$ on
$D(\E)$ if the set $\N:=\{(\phi,\dot\phi)\in 
D(W)\cap D(\hat W)\,:\, W(\phi,\dot\phi)=
\hat W(\phi,\dot\phi)\}$ is dense
in $D(\E)$. 
In the case the answer to question 1 is affirmative, two other
natural questions arise:\p 
2. Is it possible to construct such singular perturbations 
$\hat W$ without knowing $\hat A$ in advance?\p
3. Is it possible to recover the singular perturbation $\hat A$ of $A$ from
the singular perturbation $\hat W$ of $W$? 
In other words, is the following diagram commutative?
$$
\begin{CD}
W @>>> \hat W\\
@AAA   @VVV\\
A@>>> \hat A
\end{CD}
$$\p
Let us remark that in the case $\hat A$ is a strongly singular
perturbation of $A$, i.e. when the form domains of $A$ and $\hat A$ are
different, the spaces $D(\E)$ and $D(\hat \E)$ are different, so that
$W$ and $\hat W$ are defined on different
Hilbert spaces. Indeed we will answer question 2 above by looking for 
singular perturbations with 
$D(\hat\E)\simeq D(\E)\oplus (D(A)/\K)$. This results to be the
right ansatz to give affirmative answers to questions 1 and 3.\par
The framework described above can be extended by considering 
generalized abstract wave equations of the kind
$$
\ddot\phi-(C_1+C_2)\dot\phi=(A-C_1C_2)\phi\,,
$$
with both $C_1$ and $C_2$ skew-adjoint operators 
such that $A-C_1C_2$ is negative and
injective. The corresponding block matrix operator is 
$$
\overset\circ W_g=\left[\begin{matrix}
C_2&\uno\\-B^2&C_1\end{matrix}\right]
\,:\,D(B^2)\oplus D(B)\subseteq
\H_1\oplus\H_0\to \H_1\oplus\H_0\,.
$$
Then $\overset\circ W_g$ is closable, with closure $W_g$, 
as an operator on the Hilbert space 
$D(\E_C)$, the completion of $\HU\oplus\HZ$ with respect to the
scalar product
$$
\langle(\phi_1,\dot\phi_1),(\phi_2,\dot\phi_2)\rangle_{\E_C}:=
\langle B_C\phi_1,B_C\phi_2\rangle_0+
\langle\dot\phi_1,\dot\phi_2\rangle_0\,,
$$
where
$$
B_C:=(-A+C_1C_2)^{1/2}\,.
$$
Also for these generalized abstract wave equations we are able to construct 
singular perturbations $\hat W_g$ of the skew-adjoint
operator $W_g$ which reduce to the previous ones in the case
$C_1=C_2=0$. Such singular perturbation, together with their resolvents, are
defined in a relatively explicit way in terms of the original
operators $B$, $C_1$ and $C_2$. \par
The contents of the single sections are the following:\p 
{\it -- Section 2.} We review, with some variants and additions 
with respect to \cite{[K]}, \cite{[W]} and \cite{[G]} 
(and references therein) the theory of
abstract wave equations. Here we are
in particular interested (see Theorem 2.5) in computing the resolvent
of $W$, the 
skew-adjoint operator corresponding to the abstract wave
equation $\ddot\phi=-B^2\phi$, in terms of
the resolvent of $B^2$ . For such a scope the scale of Hilbert
spaces $\HB_k:=\{\phi\in\HUB\ :\ \B\phi\in D(B^{k-1})\}$, $k\ge 1$, 
is used.\p
{\it -- Section 3.} Given a continuous linear map $\tau:\HDB\to\fh$, $\fh$
an auxiliary Hilbert space, such
that, denoting by $\tau^*:\fh\to \H_{-2}$ the adjoint of the
restriction of $\tau$ to $D(B^2)$, one has Ran$(\tau^*)\cap
\H_{-1}=\emptyset$ (we are thus considering strongly singular
perturbations of $B^2$), we construct, mimiking the approch 
developed in \cite{[P1]}, skew-adjoint operators
$\hat W$ which coincide with $W$ on Ker$(\tau)\oplus\HU$. As already
mentioned, due to our hypothesis on $\tau^*$, the $\hat W$'s will be
defined on a Hilbert space larger than $\HUB\oplus\HZ$, indeed it will
a space of the kind $\HUB\oplus\HZ\oplus\fh$. Thus our strategy is the 
following: for any positive, bounded and injective self-adjoint operator 
$\Theta$ on $\fh$, at first we trivially extend $W$ to
$\HUB\oplus\HZ\oplus\fh_\Theta$ (here $\fh_\Theta$ is the Hilbert space
obtained from $\fh$ by considering the scalar product induced by
$\Theta$) by defining 
$\tilde W(\phi,\dot\phi,\zeta):=
(W(\phi,\dot\phi),0)$, which is obviously still skew-adjoint. Then we
consider the skew-symmetric operator obtained by restricting $\tilde
W$ to the kernel of the map $\tau_\Theta$, where 
$\tau_\Theta(\phi,\dot\phi,\zeta):=\tau\phi-\Theta\zeta$. 
To such a skew-symmetric operator, which depends on
$\Theta$, we apply the procedure 
given in \cite{[P1]}, thus obtaining a family of skew-adjoint
extensions parametrized by self-adjoint operators on $\fh$. Selecting
from such a family the extension corresponding to the
parametrizing operator zero, 
we obtain a skew-adjoint operator $\tilde W_\Theta$ which by
construction coincides with $\tilde W$ on the kernel of
$\tau_\Theta$ (see Theorem 3.4). 
Under the additional hypothesis that both the Hilbert spaces $\HZ$ and
$\HUB$ are contained in a common vector space (this is usually true in
the case $B$ is a (pseudo-)differential operator by considering some
space of distributions), one can then define a suitable Hilbert space
$\KUB\supset\HUB$ and a skew-adjoint operator $W_\Theta$ on
$\KUB\oplus\HZ$ such that $W_\Theta$ coincides with $W$ on the set
Ker$(\tau)\oplus\HU$ (see Theorem 3.6). 
By our hyphoteses such a set is dense in
$\HUB\oplus\HZ$ and thus $W_\Theta$ is a singular perturbation of $W$. \par
The skew-adjoint operator $W_\Theta$ permits then to define $-A_\Theta$, an injective
and positive self-adjoint operator on $\HZ$ which results to be a
singular perturbation of $-A=B^2$. The resolvent and the quadratic
form of $A_\Theta$ are also esplicitely given. Regarding the quadratic
form a variation on the Birman-Kre\u\i n-Vishik theory (see
\cite{[AS]} and references therein) is obtained. Conversely, the
skew-adjoint operator corresponding to the abstract wave equation
$\ddot\phi=A_\Theta\phi$ results to be nothing but $W_\Theta$ (there
results are summarized in Theorem 3.7). 
Thus we gave affirmative answers to questions 1-3 above.\p
{\it -- Section 4.} We construct singular perturbations of the kind obtained
in Section 3 for the
skew-adjoint operator $W_g$ corresponding now to the abstract wave equation 
$\ddot\phi-(C_1+C_2)\dot\phi=-(B^2+C_1C_2)\phi$. Here we put on the
skew-adjoint operators $C_1$ and $C_2$ conditions
which ensure that $B^2+C_1C_2$ is self-adjoint, positive and
injective. Defining $B_C:=(B^2+C_1C_2)^{1/2}$, $C:=C_1+C_2$, this case
is studied by extending the procedure of Section 3 to the 
abstract wave equation $\ddot \phi-C\dot\phi=-B^2_C\phi$ (see Theorem
4.7). The analogues of Theorems 3.4 and 3.6 corresponding to the
this more general situation are Therems 4.8 and 4.11. Here an 
hypothesis concerning both $C_1$, $C_2$ and a suitable extension
$\taub$ of the map $\tau$ must be
introduced. Such hypothesis is surely verified when $C_1$ and $C_2$
are bounded operator, whereas its validity in the unbounded case is
more subtle, as Example 3 in Section 5 shows.  \p
{\it -- Section 5.} We give some examples. In Example 1 we define
skew-adjoint operators $W_\Theta$, $\Theta$ an Hermitean injective and
positive matrix on $\C^n$, corresponding to wave equations on star-like
graphs with $n$ open ends by defining singular perturbations of the skew-adjoint operator
$W(\phi_1,\dots,\phi_n,\psi_1,\dots,\psi_n):=
(\psi_1,\dots,\psi_n,\phi''_1,\dots,\phi''_n)$, where the $\phi$'s are
defined on the half-line $(0,\infty)$ and satisfy zero Dirichlet boundary
conditions at the origin. The corresponding (according to Theorem 3.7)
negative self-adjoint operator $A_\Theta$ is of the class of
Laplacians on a star-like graphs (see \cite{[KS]} and references
therein). By a similar construction, considering also second
derivative operators on compact intervals, one could define wave equations on
more complicated graphs.\par
In Example 2 we consider the
case in which $\HZ$ is the space of square integrable functions on
$\RE^3$, $B=(-\Delta)^{1/2}$, $C_1=C_2=0$, $\fh=\C^n$ and
$\tau\phi=(\phi(y_1),\dots, \phi(y_n))$, where
$Y=\left\{y_1,\dots,y_n\right\}$ is a given discrete subset of
$\RE^3$. This gives a singular perturbations of the free wave
equations by $n$ Dirac masses placed at points $y_1,\dots, y_n$, in
the sense that the extensions constructed give a rigorous definition
and provide existence of the dynamics for wave
equations of the kind 
$$\ddot\phi=\Delta\phi+\zeta_\phi^1\delta_{y_1}+
\cdots+\zeta_\phi^n\delta_{y_n}\,,$$
where $\zeta_\phi\equiv(\zeta^1_\phi,\cdots,\zeta_\phi^n)$ is related
to the value of the continous part $\phi_0$ of $\phi$ at the points in
$Y$ by the boundary conditions 
$$\phi_0(y_i)=\sum_{1\le j\le n}\theta_{ij}\zeta^j_\phi\,,\quad
i=1,\dots,n\,.
$$ 
Such wave equations were introduced (by different methods) and
analyzed, when $n=1$, 
in \cite{[BNP]}. The corresponding singular
perturbation of the Laplacian, obtained according to Theorem 3.7 is of
the class on point perturbation of the Laplacian (see \cite{[AGHH]}
and references therein). The above situation can be generalized by 
taking as $\tau$ the evaluation map along a $d$-set (i.e. a
$d$-dimensional Lipschitz
submanifold if $d$ is an integer or a self-similar fractal in the
noninteger case), with  and proceeding similarly to the
examples appearing in \cite{[P1]}-\cite{[P3]}, 
thus obtaining perturbations of the
free wave equation supported on null sets. Here the extension
paramenter is a self-adjoint operator on some
fractional order Sobolev space on the $d$-set. \par 
A wave equation of the kind
$\ddot\phi=\Delta\phi+4\pi e\,M\zeta_\phi\delta_0$ was used 
to give a rigorous description of classical and quantum
electrodynamic in dipole (or linear) approximation and without ultraviolet cut-off 
(see \cite {[NP]} and \cite{[BNP2]}). Here $\phi$ is
$\RE^3$-valued and plays the
role of the electromagnetic potential in the Coulomb gauge 
(thus div $\phi=0$), $M$ is the projector onto the divergenceless
fields and $e$ is the electric charge (the velocity
of light being set to be equal to one). In this case one must
modify the above boundary condition (here $Y=\left\{0\right\}$), 
considering the (no more linear but affine) one
given by
$$
\phi_0(0)=-\frac{m}{e}\,\zeta_\phi+\frac{1}{e}\, p\,,
$$ 
where $p$ is an arbitrary vector in $\RE^3$ and $m$ is the mass of the
particle. In this framework $\zeta_\phi\in\RE^3$ can be identified with the
particle velocity $v$, so that the particle dynamics is given by the
evolution of the field
singularity. With this identification 
the above boundary condition is nothing else that the 
usual (linearized and regularized) relation between velocity and momentum 
(represented by the vector
$p$) in the presence on an electromagnetic field, i.e. $p=mv+e\,\phi_0(0)$.\par
This approach suggests that the study of singular perturbations
of the wave equation $\ddot\phi=\Delta\phi$ can produce an useful
framework for a rigorous treatment of classical
electrodynamics of point particles and for quantum electrodynamics in 
the ultraviolet limit. Indeed this was the original motivation of the
paper. In order to remove the limitation given by the
dipole approximation assumed in \cite{[NP]} and \cite{[BNP2]}, one is lead to study the
singular perturbations, supported at the origin, 
of the wave equation 
\begin{align*}
&\dot\phi=v\cdot\nabla\phi+\psi\\
&\dot\psi=v\cdot\nabla\psi+\Delta\phi\,,\\
\end{align*}
were $v$ is a given vector in $\RE^3$ with $|v|<1$. This is suggested by starting with the Maxwell-Lorentz
system, by re-writing it in a reference frame co-moving with the particle
and then by performing the reduction allowed by the conservation of
the total (particle + field) momentum. We refer
to the digression given at the end of Section 5 for a more detailed 
discussion. Thus in the successive example in
Section 5 (Example 3), we modify the situation considered in
Example 2 (in the case $Y=\left\{0\right\}$) by taking
$C_1=C_2=v\cdot\nabla$, with $v\in\RE^3$, $|v|<1$. In this case the
regular part $\phi_0$
of $\phi$'s in the
proper operator domain is no more continuous (when $v\not=0$) 
and the elavuation map
of Example 2 has to be extended to
$\taub$, where $\taub\phi_0$ is defined by the 
limit $R\downarrow 0$ of the average $\langle\phi_0\rangle_R$
of $\phi_0$ over the sphere of radius $R$. It is here proven the 
such a limit exists for the functions in the operator domain of the 
extensions. This produces 
a rigorous definition and existence of the dynamics for the wave equation
\begin{align*}
&\dot\phi=v\cdot\nabla\phi+\psi\\
&\dot\psi=v\cdot\nabla\psi+\Delta\phi+\zeta_\phi\delta_0\,,\\
\end{align*}
where now the $\zeta_\phi$'s are related to the regular part $\phi_0$
of the $\phi$'s by the boundary condition
$$
\langle\phi_0\rangle:=\lim_{R\downarrow 0}\,\langle\phi_0\rangle_R=\theta\zeta_\phi\,.
$$
Once the proper domain of definition for the fields $\phi$ and $\psi$ 
is determined by 
this linear analysis, a nonlinear
operator, candidate to describe the classical electrodynamics of a
point 
particle,
can be obtained by considering the nonlinear wave equation 
\begin{align*}
&\dot\phi=v\cdot\nabla\phi+\psi\\
&\dot\psi=v\cdot\nabla\psi+\Delta\phi+4\pi e\,Mv\delta_0\,,\\
\end{align*}
where $v$, again representing the particle velocity, 
is no more a given vector but is
related to the regular parts $\phi_0$ and $\psi_0$ of the fields
$\phi$ and $\psi$ by the nonlinear boundary condition
$$
\langle\phi_0\rangle-\frac{1}{4\pi
e}\,\langle\psi_0,\nabla\phi_0\rangle=-\frac{m}{e}\,
\frac{v}{\sqrt{1-|v|^2}}+\frac{1}{e}\,\Pi\,.
$$ 
The (conserved) total momentum $\Pi$ of the particle-field system
is defined, in terms of the particle momentum $p$, by
$\Pi:=p-\frac{1}{4\pi}\,\langle\psi,\nabla\phi\rangle$. Thus 
the above boundary condition corresponds to 
the (regularized) velocity-momentum relation for a (relativistic)
particle in the presence of an electromagnetic field, i.e. 
$p=\frac{m v}{\sqrt{1-|v|^2}}+e\,\langle\phi_0\rangle$. Again we refer
to the digression at the end of Section 5 for more details.
\p
{\it -- Appendix.} We give a compact rewiev of the approach to singular
perturbations of self-adjoint operators developed in \cite{[P1]}
adapted to our present (skew-adjoint) situation. In particular, with
reference to the notations in \cite{[P1]}, we make here
a particular choice of the operator $\Gamma$ which correspond, 
in the case treated in Section 3 here, to a weakly singular
perturbation. Thus a strongly singular perturbation $\hat A$ of $A$ gives rise
to a weakly singular perturbation $\hat W$ of $W$. This could be used
to study the scattering theory for strongly singular perturbations of
$A$ in terms of weakly singular perturbations. Indeed, by Birman-Kato
invariance principle, the M\"oller operators $\Omega_\pm(\hat W,W)$
and $\Omega_\pm(\hat A,A)$ are unitarily equivalent. 
As regard the parametrizing operator,
as we already said above, we pick up here, in the family of
skew-adjoint extensions given by the general scheme in \cite{[P1]},
the extension corresponding to the zero operator.  
\def\CD{\bar C_2}

\section{abstract wave equations}
Let $B: D(B)\subseteq\HZ\to \HZ$ be a self-adjoint operator on the Hilbert space $\HZ$ such that 
Ker$(B)=\left\{0\right\}$. Let us denote by $\H_k$, $k>0$, the scale
of Hilbert spaces
given by the domain of $B^k$ with the scalar product 
$\langle\cdot,\cdot\rangle_k$ leading to the
graph norm, i.e.
$$
\langle\phi_1,\phi_2\rangle_k:=
\langle B^k\phi_1,B^k\phi_2\rangle_0+\langle \phi_1,\phi_2\rangle_0\,.
$$
Here $\langle\cdot,\cdot\rangle_0$ denotes the scalar product in
$\HZ$. We will use the symbol $\|\cdot\|_0$ to indicate the
corresponding norm.\par
We then define the Hilbert space $\HUB$ by completing 
the pre-Hilbert space $D(B)$ endowed with the scalar product
$$
[\phi_1,\phi_2]_1:=
\langle B\phi_1,B\phi_2\rangle_0\,.
$$
We define $\B\in \BO(\HUB,\HZ)$ as the closed bounded extension of the densely
defined linear operator
$$
B:\HU\subseteq\HUB\to\H\,.
$$ 
Here and below by $\BO(X,Y)$ we mean the space of bounded, everywhere
defined, linear operators on the Banach space $X$ to the Banach space
$Y$; for brevity we put $\BO(X)\equiv\BO(X,X)$.\par
Since $B$ is self-adjoint one has 
$$
\ran(B)^\perp=\ker(B)\,,
$$
so that, $B$ being injective, $\ran(B)$ is dense in $\HZ$. Therefore we
can define $\B^{-1}\in\BO(\HZ,\HUB)$ as the closed bounded extension
of the densely defined linear operator
$$
B^{-1}:\ran(B)\subseteq\HZ\to\HUB\,.
$$
One can then verify that $\B$ is boundedly invertible with inverse
given by $\B^{-1}$. \par
Given $\B$ we introduce the scale of spaces $\HB_k$, $k\ge 1$, 
defined by
$$
\HB_k:=\left\{\phi\in\HUB\ :\ \B\phi\in\H_{k-1}\right\}\,.
$$
Obviously $\H_k\subseteq \HB_k$. 
\begin{lemma} 
$$\HB_k=\H_k+\HB_{k+1}\,.$$
\end{lemma}
\begin{proof} 
The thesis follow from 
\begin{align*}
\HB_{2k}=&\B^{-1}(B+i)^{-1}(B^{2(k-1)}+1)^{-1}(\H_0)\,,\\
\H_{2k}=&(B^2+1)^{-1}(B^{2(k-1)}+1)^{-1}(\H_0)\,,
\end{align*}
\begin{align*}
\HB_{2k+1}=&\B^{-1}(B^{2k}+1)^{-1}(\H_0)\,,\\
\H_{2k+1}=&(B+i)^{-1}(B^{2k}+1)^{-1}(\H_0)\,,
\end{align*}
and from the identities
$$
\B^{-1}=(B+i)^{-1}+i\B^{-1}(B+i)^{-1}\,,
$$
$$
\B^{-1}(B+i)^{-1}=(B^2+1)^{-1}-i\B^{-1}(B^2+1)^{-1}\,.
$$
\end{proof}
\begin{lemma} The set $\HB_k$ endowed with the scalar product 
$$
[\phi_1,\phi_2]_k:=\langle\B\phi_1,\B\phi_2\rangle_{k-1}
$$
is a Hilbert space.
\end{lemma} 
\begin{proof} 
Let $\phi_n$, $n\ge 1$, be a Cauchy sequence
in $\HB_k$. Then $\phi_n$, $n\ge 1$, is Cauchy in
$\HUB$ and $\B\phi_n$, $n\ge 1$, is Cauchy in $\HB^{k-1}$. 
Thus $\B\phi_n\to\B\phi$ and 
$B^{k-1} \B\phi_n\to\psi$ in $\HZ$. 
Since $B^{k-1}$ is closed, $\B\phi\in\H_{k-1}$, hence 
$\phi\in\HB_{k}$, and $\psi=B^{k-1}\B\phi$. 
\end{proof}
\begin{remark} 
The previous lemma shows that $\HB_k$ could be alternatively defined as the
completion of pre-Hilbert space $D(B^k)$ endowed with the scalar product
$$
[\phi_1,\phi_2]_k:=\langle B\phi_1, B\phi_2\rangle_{k-1}\,.
$$
Thus $\H_k$ is dense in $\HB_k$. 
\end{remark}
We now define
$$
\A:\HDB\to\HZ\,,\qquad \A:=-B\B\,.
$$ 
\begin{remark} By the previous remark $\A\in\BO(\HDB,\HZ)$ 
could be alternatively defined as the
closed bounded extension of the densely defined linear operator 
$A:=-B^2:\HD\subseteq\HDB\to\HZ$.
\end{remark}
We put, for any real $\lambda\not=0$, 
$$
R_0(\lambda):=(B^2+\lambda^2)^{-1}\,,\qquad R_0(\lambda)\in\BO(\HZ,\HD) 
$$
and then define $\R\in\BO(\HUB,\HB_3)$ as the closed bounded extension of
$$
R_0(\lambda):\HU\subseteq\HUB\to\HB_3\,.
$$
The linear operator $\R$ 
satisfies the relations
\begin{equation}
-\A\R+\lambda^2\R=\uno_{\HUB}\,,
\end{equation}
\begin{equation}
-R_0(\lambda)\A+\lambda^2\R=\uno_{\HDB}\,,
\end{equation}
On the Hilbert space $\HUB\oplus\HZ$ with scalar product given by 
$$
\langle\langle\,(\phi_1,\psi_1),(\phi_2,\psi_2)\,\rangle\rangle:=
\langle\B\phi_1,\B\phi_2\rangle_0+\langle\psi_1,\psi_2\rangle_0\,.
$$
we define the linear operator
$$
W:\HDB\oplus\HU\subseteq\HUB\oplus\HZ\to \HUB\oplus\HZ\,,\qquad 
W(\phi,\psi):=(\psi,\A\phi)\,.
$$
\begin{theorem} The linear operator $W$ is skew-adjoint and its
resolvent is given by
$$
(-W+\lambda)^{-1}(\phi,\psi)
=(\lambda\R\phi+R_0(\lambda)\psi,-\phi+\lambda^2\R\phi+\lambda R_0(\lambda)\psi)\,.
$$
\end{theorem}
\begin{proof} The skew-symmetry of $W$ immediately follows from the
definition of the scalar product 
$\langle\langle\cdot,\cdot\rangle\rangle$. The fact that 
$(-W+\lambda)^{-1}$ as defined above is the inverse of $-W+\lambda$ is a
matter of algebraic computations given the definition of $R_0(\lambda)$,
$\R$ and (2.1), (2.2). The proof is then concluded by recalling that $W$ is
skew-adjoint (equivalently $iW$ is self-adjoint) if and only if it is 
skew-symmetric and Ran$(W\pm
\lambda)=\HUB\oplus\HZ$ for some real $\lambda\not= 0$. 
\end{proof} 
\begin{remark} Note that $\HDB\oplus\HU=\ran(-W+\lambda)^{-1}$ gives a 
decomposition compatible with the one given
by lemma 1.1, i.e. $\HDB=\HD+\HB_3$ and $\HU=\HUB+\HD+\HB_3$.
\end{remark}
\begin{remark} 
Note that the norm on $\HDB$ induced by the graph norm of $W$ coincides with
the one given by the scalar product $[\cdot,\cdot]_2$. Hence the
domain of $W$ is the direct sum of the Hilbert spaces $\HDB$ and $\HU$
as written above.
\end{remark}

\section{Singular perturbations of abstract wave equations}
On the Hilbert space $\fh$ with scalar product 
$\langle\cdot,\cdot\rangle_\fh$ and norm $\|\cdot\|_\fh$, we consider
a bounded, positive and injective self-adjoint operator $\Theta$.
Then we denote by $\fh_\Theta$ the
Hilbert space given by $\fh$ endowed with the scalar product
$$
\langle\zeta_1,\zeta_2\rangle_\Theta:=
\langle\Theta \zeta_1,\zeta_2\rangle_\fh\,.
$$ 
The corresponding norm will be indicated by $\|\cdot\|_\Theta$. \par
By Theorem 2.5, on Hilbert space $\HUB\oplus\HZ\oplus\fht$ with
scalar product 
$$
\langle\langle\,(\phi_1,\psi_1,\zeta_1),(\phi_2,\psi_2,\zeta_2)\,\rangle\rangle_{\Theta}:=
\langle\B\phi_1,\B\phi_2\rangle_0+\langle\psi_1,\psi_2\rangle_0+
\langle \zeta_1,\zeta_2\rangle_\Theta\,,
$$
the linear operator
$$
\tilde W:\HDB\oplus\HU\oplus\fht\subseteq\HUB\oplus\HZ\oplus\fht
\to \HUB\oplus\HZ\oplus\fht\,,
$$
\begin{equation*}
\tilde W(\phi,\psi,\zeta):=(W(\psi,\phi),0)
\end{equation*}
is skew-adjoint and 
\begin{equation}
(-\tilde W+\lambda)^{-1}(\phi,\psi,\zeta)
=((-W+\lambda)^{-1}(\phi,\psi),\lambda^{-1}\zeta)\,.
\end{equation}
Given $\tau\in\BO(\HDB,\fh)$, we define $\taut\in
\BO(\HDB\oplus\HZ\oplus\fht,\fh)$ by
$$
\taut:\HDB\oplus\HZ\oplus\fht\to\fh\,,\qquad
\taut(\phi,\psi,\zeta):=
\tau\phi-\Theta \zeta\,.
$$
The action of $\taut$ satisfies A.1 (see the appendix). 
Now we suppose that it also
satisfies A.2, i.e. we suppose \p
H3.0)
$$
\ran(\tau_\Theta)=\fh\,.
$$
Of course H3.0 holds true if $\tau$ itself is surjective. Another 
possibility is  
$$
\forall\,\zeta\in\fh\,,\quad
\|\Theta\zeta\|_\fh\ge c\,\|\zeta\|_\fh\,,\quad c>0\,,
$$ 
which is equivalent to $\ran(\Theta)=\fh$.\par
Now we define $\Gb\in\BO(\HZ,\fh)$ and $\G\in\BO(\fh,\HZ)$ by
$$
\Gb:=\tau R_0(\lambda)\,,\qquad \G:=\Gb^*\,.
$$ 
We also define $\GBb\in\BO(\HUB,\fh)$ and $\GB\in\BO(\fh,\HUB)$ by
$$
\GBb:=\tau\R\,,\qquad\GB:=\GBb^*\,.
$$  
Obviously $\GBb=\Gb$ on $\HU$.
\begin{lemma} 
$$\ran(\G)\subseteq\ran(B)\quad\text{ and}\quad 
\GB=\B^{-1} B^{-1}\G\in\BO(\fh,\HDB)\,.$$
\end{lemma}
\begin{proof} By the definitions of $\GB$ and $\G$ one has, for any
$\zeta\in\fh$ and for any $\psi\in\HU$,
$$
\langle B\psi,\B\GB \zeta\rangle_0=\langle\Gb\psi,\zeta\rangle_\fh=
\langle \psi,\G\zeta\rangle_0\,.
$$ 
Being $B$ self-adjoint with domain $\HU$, the above relation shows
that $\B\GB \zeta\in\HU$, hence $\GB \zeta\in\HDB$,  
\begin{equation*}
B\B\GB=\G\,,
\end{equation*}
$$
\|\B\GB \zeta\|_0^2+\|B\B\GB \zeta\|^2_0=
\|\B\GB \zeta\|_0^2+\|\G\zeta\|^2_0\,.
$$
\end{proof}
Defining 
$\GTb\in\BO(\HUB\oplus\HZ\oplus\fht,\fh)$ by
\begin{equation}
\GTb(\phi,\psi,\zeta):=
\taut(-\tilde W+\lambda)^{-1}(\phi,\psi,\zeta)
=\lambda\GBb\phi+\Gb\psi-\lambda^{-1}\Theta \zeta\,
\end{equation}
and
$\GT\in\BO(\fh,\HUB\oplus\HZ\oplus\fht)$ by
\begin{equation}
\GT \zeta:=-\,\breve G_\Theta(-\lambda)^*\zeta
=(\lambda\GB \zeta,-\,\G \zeta,-\lambda^{-1}\zeta)\,,
\end{equation}
one has that, by the previous lemma,
$\GT\in\BO(\fh,\HDB\oplus\HZ\oplus\fht)$. Thus A.4 is satisfied,
\begin{align}
\GMT:=&-\,\taut\GT=-\lambda\tau\GB-\frac{1}{\lambda}\,\Theta\\
=&-\lambda\tau\B^{-1} B^{-1} G(\lambda)
-\frac{1}{\lambda}\,\Theta
\end{align}
is well-defined and $\GMT\in\BO(\fh)$. Let us now show that 
A.5 is satisfied:
\begin{lemma}
$$\Gamma_\Theta(\lambda)^*=-\Gamma_\Theta(-\lambda)\,.$$
\end{lemma}
\begin{proof} By \cite{[P1]}, Lemma 2.1, 
$$
(\lambda^2-\epsilon^2)\,R_0(\epsilon)\,G(\lambda)=G(\epsilon)-G(\lambda)\,.
$$
Since $\ran(\G)\subseteq\ran(B)$, $R_0(\epsilon)\,G(\lambda)$ strongly
converges in $\BO(\fh,\HDB)$, as $\epsilon\downarrow 0$, to 
$\B^{-1} B^{-1} G(\lambda)$ when $B^2R_0(\epsilon)$ strongly converges to
the identity operator on $\HZ$. Since $B^2$ is injective this follows
proceeding as in \cite{[P2]}, Section 3. Therefore one has that     
\begin{align*}
\Gamma_\Theta(\lambda)=&\text{\rm s-}\lim_{\epsilon\downarrow 0}\, 
-\frac{1}{\lambda}\,\left(\Theta+\tau(G(\epsilon)-G(\lambda))\right)\\
=&\text{\rm s-}\lim_{\epsilon\downarrow 0}\, 
-\frac{1}{\lambda}\,\left(\Theta+(\lambda^2-\epsilon^2)\,
\breve G(\epsilon)G(\lambda)\right)\\
=&\text{\rm s-}\lim_{\epsilon\downarrow 0}\, 
-\frac{1}{\lambda}\,\left(\Theta+(\lambda^2-\epsilon^2)\,
\breve G(\lambda)G(\epsilon)\right)
\,.
\end{align*}
The proof is the concluded by observing that
$\tau(G(\epsilon))-G(\lambda))$ is symmetric (see \cite{[P1]}, Lemma
2.2. Also see \cite{[P2]}, Lemma 3).
\end{proof}
\begin{remark} By the same methods used in the above proof (i.e. using
the fact that $\ran(\G)\subseteq\ran(B)$), all the results contained in
\cite{[P2]} can be extended to the case in which $\tau\in\BO(\HDB,\fh)$, thus
allowing for the treatment of singular perturbations of convolution
operators also in lower dimensions (in \cite{[P2]} the examples were given in
$\RE^d$ with $d\ge 4$).
\end{remark}
Denote by $\H_{-k}$, $k\ge 0$, the completion of $\HZ$ with respect to 
the scalar product 
$$
\langle\phi_1,\phi_2\rangle_{-k}
:=\langle(B^{2k}+1)^{-1/2}\phi_1,(B^{2k}+1)^{-1/2}\phi_2\rangle_0\,.
$$
Of course $\H_{-k}\subseteq\H_{-(k+1)}$. Since
$\tau\in\BO(\HD,\fh)$ we define $\tau^*\in\BO(\fh,\H_{-2})$ by 
$$
\langle(B^{4}+1)^{-1/2}\tau^*\zeta,(B^{4}+1)^{1/2}\phi\rangle_0=
\langle\zeta,\tau\phi\rangle_\fh\,,\quad \zeta\in\fh\,,\,\phi\in\HD\,.
$$
Now we suppose that \p
H3.1)
$$
\ran(\tau^*)\cap\H_{-1}=\{0\}\,.
$$
This, using the definition of $\G$, is equivalent to
$$
\ran(\G)\cap\H_1=\left\{0\right\}\,,
$$
so that A.3 is satisfied, i.e.
$$
\ran(\GT)\cap D(\tilde W)=\left\{0\right\}\,.
$$
By Theorem 6.2 we can define a skew-adjoint extension of 
the skew-symmetric operator 
given by the restriction of $\tilde W$ to the dense set
$$\N_\Theta:=\left\{(\phi,\psi,\zeta)\in \HDB\oplus\HU\oplus\fht\ :\
\tau\phi=\Theta \zeta\right\}\,:
$$
\begin{theorem} Suppose that \text{\rm H3.0} and \text{\rm H3.1} hold
true. Let 
\begin{align*}
&D(\tilde W_\Theta):=\left\{(\phi_0,\psi,\zeta_\phi)\in
\HUB\oplus\HZ\oplus\fht\ :\quad \phi_0\in\HDB\,,\right.\\
&\left.\psi=\psi_\lambda+\G \zeta_\psi\,,
\quad\psi_\lambda\in\HU\,,\ 
\zeta_\psi\in\fh\,,\quad \Theta \zeta_\phi=\tau\phi_0\right\}\,.
\end{align*}
Then
$$
\tilde W_\Theta:D(\tilde W_\Theta)\subseteq 
\HUB\oplus\HZ\oplus\fht\to
\HUB\oplus\HZ\oplus\fht\,,\qquad$$
$$\tilde W_\Theta(\phi_0,\psi,\zeta_\phi):=(\psi_0,\,\A\phi_0,\,\zeta_\psi)\,,
$$
is a skew-adjoint extension of the restriction of $\tilde W$ to 
the dense set $\N_\Theta$. Here $\psi_0\in\HUB$, defined by 
$$\psi_0:=\psi_\lambda-\lambda^2\B^{-1} B^{-1} \G\,
\zeta_\psi\,,
$$ does not depend on $\lambda$. 
The resolvent of $\tilde W_\Theta$ is given by
$$
(-\tilde W_\Theta+\lambda)^{-1}=(-\tilde W+\lambda)^{-1}+\GT\GMT^{-1}\GTb\,,
$$
where the bounded linear operators $(-\tilde W+\lambda)^{-1}$, $\GTb$, $\GT$, 
$\GMT^{-1}$ have been defined in \text{\rm (3.1)-(3.4)} respectively.
\end{theorem}
\begin{proof} By Theorem 6.4 we known that 
$(-\tilde W+\lambda)^{-1}+\GT\GMT^{-1}\GTb$ is the resolvent
of a skew-adjoint extension $\hat W_\Theta$ of 
the restriction of $\tilde W$ to 
the dense set $\N_\Theta$. Therefore $(\hat \phi_0,\hat
\psi,\hat \zeta_\phi)\in D(\hat W_\Theta)$
if and only if
\begin{align*}
\hat\phi_0=&\,\phi_\lambda+\lambda\GB\GMT^{-1}(\tau\phi_\lambda-\Theta
\zeta_\lambda)\,,\qquad \phi_\lambda\in\HDB\,,\\
\hat\psi=&\,\psi_\lambda-\G\GMT^{-1}(\tau\phi_\lambda-\Theta
\zeta_\lambda)\,,\qquad \psi_\lambda\in\HU\,,\\
\hat \zeta_\phi=&\,\zeta_\lambda-\frac{1}{\lambda}\,\GMT^{-1}(\tau\phi_\lambda-\Theta
\zeta_\lambda)\,,\qquad \zeta_\lambda\in\fh\,.
\end{align*} 
Let us now show that $D(\hat W_\Theta)=D(\tilde W_\Theta)$. \par
Since $\ran(\GB)\subseteq\HDB$, so that $\hat\phi_0\in\HDB$, and 
\begin{align*}
&\taut((-\tilde W+\lambda)^{-1}+\GT\GMT^{-1}\GTb)\\
&=
\GTb-\GMT\GMT^{-1}\GTb=0, 
\end{align*}
so that $\tau\hat\phi_0=\Theta \hat \zeta_\phi$, we have $D(\hat
W_\Theta)\subseteq D(\tilde W_\Theta)$. Let us now prove the reverse
inclusion. Given $(\phi_0,\psi,\zeta_\phi)\in D(\tilde W_\Theta)$ let us define 
\begin{align*}
\phi_\lambda:=&\,\phi_0+\lambda\GB\,\zeta_\psi\,,\\
\zeta_\lambda:=&\,\zeta-\frac{1}{\lambda}\,\zeta_\psi\,.
\end{align*}
Then
$$
\tau\phi_0
=\tau\phi_\lambda-\lambda\tau\GB\,\zeta_\psi
=\Theta \zeta=\Theta\left(\zeta_\lambda+\frac{1}{\lambda}\,\zeta_\psi\right)
$$
implies
$$
\tau\phi_\lambda-\Theta \zeta_\lambda=\left(
\lambda\tau\GB+\frac{1}{\lambda}\,\Theta\right)\zeta_\psi\,,
$$
i.e.
$$
\zeta_\psi=-\,\GMT^{-1}(\tau\phi_\lambda-\Theta \zeta_\lambda)\,.
$$
Thus $D(\tilde W_\Theta)\subseteq D(\hat W_\Theta)$. Now we have
\begin{align*}
\hat W_\Theta(\phi_0,\psi,\zeta_\phi)=&\tilde
W(\phi_\lambda,\psi_\lambda,\zeta_\lambda)+\lambda
(\phi_0-\phi_\lambda,\psi-\psi_\lambda,\zeta_\phi-\zeta_\lambda)\\
=&(\psi_\lambda-\lambda^2\GB\,\zeta_\psi,\A\phi_\lambda+\lambda\G
\zeta_\psi,\zeta_\psi)\\
=&(\psi_0,\A(\phi_\lambda-\lambda\GB\,\zeta_\psi),\zeta_\psi)\\
=&(\psi_0,\A\phi_0,\zeta_\psi)\\
=&\tilde W_\Theta(\phi_0,\psi,\zeta_\phi)\,.
\end{align*}
$\psi_0$ does not depend on $\lambda$ since the definition of 
$\tilde W_\Theta$ is $\lambda$-independent.
\end{proof}
Let us now suppose that \smallskip\p
H3.2) \text{\it $\quad$both $\H_0$ and $\HUB$ 
are contained in a given vector space
$\V$.}
\smallskip\p
Thus we can define 
$$
G:\fh\to\V\,,\qquad G:=\G+\lambda^2\GB\,.
$$
\begin{lemma} The definition of $G$ is $\lambda$-independent. Moreover 
$$
\ran(G)\cap\HUB=\left\{0\right\}\,.
$$
\end{lemma}
\begin{proof} By first resolvent identity one has (see \cite{[P1]}, Lemma 2.1)
$$
(\lambda^2-\mu^2)\,R_0(\mu)\,G(\lambda)=G(\mu)-G(\lambda)\,,
$$
i.e.
$$
\lambda^2G(\lambda)-\mu^2G(\mu)=B^2(G(\mu)-G(\lambda))\,.
$$
This implies, by Lemma 3.1,
$$
\G+\lambda^2\GB=G(\mu)+\mu^2\bar G(\mu)\,.
$$
Suppose there exists $\zeta\in\fh$ such that
$$
\G\zeta+\lambda^2\GB\zeta=\phi\in\HUB\,.
$$
Then $\G\zeta\in\HU$ and so, by H3.1, $\G\zeta=0$. By Lemma 3.1
$\GB\zeta=0$ and the proof is done.
\end{proof}
By the previous lemma the following spaces are well-defined: 
$$
\KUB:=\left\{\phi\in\V\ :\ \phi=\phi_0+G\zeta_\phi\,,\ \phi_0\in\HUB\,,\
\zeta_\phi\in\fh\right\}\,,
$$
$$
\KDB:=\left\{\phi\in\V\ :\ \phi=\phi_0+G\zeta_\phi\,,\ \phi_0\in\HDB\,,\
\zeta_\phi\in\fh\right\}\,,
$$
$$
\KU:=\KUB\cap\HZ\,.
$$
Moreover the map 
$$
U:\HUB\oplus\HZ\oplus\fh_\Theta\to\KUB\oplus\HZ\,,\qquad 
U(\phi_0,\psi,\zeta_\phi):=(\phi_0+G\zeta_\phi,\psi)
$$
is injective and surjective and thus is unitary once 
we make $\KUB$ a Hilbert space by defining the scalar product
$$
\langle\phi,\varphi\rangle_{\KUB}:=
\langle\B\phi_0,\B\varphi_0\rangle_0
+\langle\zeta_\phi,\zeta_\varphi\rangle_\Theta\,.
$$
Thus we can state the following:
\begin{theorem} Suppose that \text{\rm H3.0}, \text{\rm H3.1} and 
\text{\rm H3.2} hold true.
Then the linear operator 
$$
W_\Theta:D(W_\Theta)\subseteq\KUB\oplus\HZ\to\KUB\oplus\HZ\,,
$$
$$
D(W_\Theta)=\left\{(\phi,\psi)\in\KDB\oplus\KU\ :\ 
\Theta\,\zeta_\phi=\tau\phi_0\right\}\,,
$$
$$
W_\Theta(\phi,\psi):=U \tilde W_\Theta U^*(\phi,\psi)=
(\psi,\A\phi_0)\,.
$$
is skew-adjoint. It coincides with 
$$W:\HDB\oplus\HU\subseteq\HUB\oplus\HZ\to\HUB\oplus\HZ\,,
\quad W(\phi,\psi)=(\psi,\A\phi)$$ 
on the dense set 
$$D(W)\cap D(W_\Theta)=\left\{\phi\in\HDB\ :\
\tau\phi=0\right\}\oplus\HU\,. $$
\end{theorem}
Once we obtained $W_\Theta$ we can define the linear operator
$A_\Theta$ on $\HZ$ by
$$
D(A_\Theta):=\left\{\phi\in\KDB\cap\HZ\ :\ \Theta\zeta_\phi=\tau\phi_0\right\}\,,
$$
$$
A_\Theta :D(A_\Theta)\subseteq\HZ\to\HZ\,,\qquad A_\Theta\phi:=
P_2W_\Theta I_1\phi\equiv\bar A\phi_0\,,
$$
where
$$
P_2:\KUB\oplus\HZ\to\HZ\qquad P_2(\phi,\psi):=\psi\,,
$$ 
and  
$$
I_1:\KDB\cap\HZ\to\KDB\oplus\KU\,,\qquad I_1\phi:=(\phi,0)\,.
$$ 
We have the following
\begin{theorem} 1. $A_\Theta$ is a negative and injective self-adjoint
operator which coincides with $A$ on the set $\ker (\tau)$. Its
resolvent is given by
$$
(-A_\Theta+\lambda^2)^{-1}=R_0(\lambda)+
\G(\Theta+\lambda^2\tau\B^{-1}B^{-1}\G)^{-1}\Gb\,.
$$
The positive quadratic form $\Q$ corresponding to $-A_\Theta$ is
$$
\Q:\K_1\subseteq\HZ\to\RE\,,\quad \Q(\phi)=\|\B\phi_0\|^2_0+
\|\zeta_\phi\|^2_{\Theta}\,.
$$
2. The skew-adjoint operator corresponding to the abstract wave
   equation $\ddot\phi=A_\Theta\phi$ is the skew-adjoint operator
   $W_\Theta$ defined in the previous theorem. 
\end{theorem} 
\begin{proof} 1. Let us define
$$
R_\Theta(\lambda):=
R_0(\lambda)+
\G(-\lambda\Gamma_\Theta(\lambda))^{-1}\Gb\,.
$$
By the proof of Lemma 3.2 and \cite{[P1]}, Lemma 2.1, 
\begin{align*}
&-\lambda\Gamma_\Theta(\lambda)-(-\mu\Gamma_\Theta(\mu))\\
=&\text{\rm s-}\lim_{\epsilon\downarrow 0}\, 
(\lambda^2-\epsilon^2)\,
\breve G(\lambda)G(\epsilon)-(\mu^2-\epsilon^2)\,
\breve G(\mu)G(\epsilon)\\
=&\tau(G(\mu)-G(\lambda))=(\lambda^2-\mu^2)\breve G(\mu)\Gb\,.
\end{align*}
We already know that $-\lambda\Gamma_\Theta(\lambda)$ is boundedly
invertible and, by (3.5) and Lemma 3.2, $(-\lambda\Gamma_\Theta(\lambda))^*=
-\lambda\Gamma_\Theta(\lambda)$. Therefore, by
\cite{[P1]}, Proposition 2.1, $R_\Theta(\lambda)$ is the resolvent of
a self-adjoint operator $\tilde A_\Theta$, coinciding with $A$ on
$\ker(\tau)$, defined by
$$
D(\tilde A_\Theta):=\left\{\phi\in\HZ\,:\,\phi=\phi_\lambda
+\G(-\lambda\Gamma_\Theta(\lambda))^{-1}\tau\phi_\lambda\right\}
$$ 
$$
(-\tilde A_\Theta+\lambda^2)\phi:=(-A+\lambda^2)\phi_\lambda\,.
$$ 
One then proves that $\tilde A_\Theta\equiv A_\Theta$ proceeding
 exactly as in the proof of \cite{[P2]}, Theorem 5.
\par
Since $A$ is injective, $A_\Theta\phi=0$ implies 
$\phi_0=0$ and thus $\zeta_\phi=0$, i.e. $\phi=0$. \par
By the proof of Lemma 3.1 one has
$$
\langle B\phi_0,\B\GB \zeta\rangle_0=\langle\GBb\phi_0,\zeta\rangle_\fh
$$
and, by (2.2), 
$$
\langle B\B\phi_0,\G\zeta_\phi\rangle_0
+\lambda^2[\phi_0,\GB\zeta_\phi]_1
=\langle\tau\phi_0,\zeta_\phi\rangle_\fh\,.
$$
Thus, using the definition of $G$ and the two different decompositions 
of $\phi\in D(A_\Theta)$ given by 
$$
\phi=\phi_0+G\zeta_\phi=\phi_\lambda+\G\zeta_\phi\,,
$$
one obtains
\begin{align*}
&\langle -A_\Theta\phi,\phi\rangle_0=
\langle -\A\phi_0,\phi_\lambda\rangle_0+
\langle -\A\phi_0,\G\zeta_\phi\rangle_0\\
=&\langle \B\phi_0,\B\phi_\lambda\rangle_0+
\langle \B B\phi_0,\G\zeta_\phi\rangle_0\\
=&\langle \B\phi_0,\B\phi_0\rangle_0+\lambda^2
\langle \B\phi_0,\B\G\zeta_\phi\rangle_0
-\lambda^2\langle \B\phi_0,\B\G\zeta_\phi\rangle_0
+\langle\tau\phi_0,\zeta_\phi\rangle_\fh\\
=&\langle \B\phi_0,\B\phi_0\rangle_0
+\langle\zeta_\phi,\zeta_\phi\rangle_\Theta\,.
\end{align*}
Thus $A_\Theta$ is negative. Since $\K_1$ is obviously complete with
respect to the norm 
$$
\|\phi\|^2_{\K_1}:=\|\B\phi_0\|^2_0
+\|\zeta_\phi\|^2_\Theta
+\|\phi\|^2_0\,,
$$
the closed and positive quadratic form $\Q$ is the one 
associated to $-A_\Theta$. \p
2. Since the completion of $\KU$ with respect to the scalar
product $$[\phi,\varphi]_1:=
\langle B\phi_0,B\varphi_0\rangle+\langle\zeta_\phi,\zeta_\varphi\rangle_\Theta
$$
is $\KUB$ and the completion of $D(A_\Theta)$ 
with respect to the scalar
product $$[\phi,\varphi]_2:=
\langle B\phi_0,B\varphi_0\rangle
+\langle\zeta_\phi,\zeta_\varphi\rangle_\Theta+
\langle A\phi_0,A\varphi_0\rangle
$$
is $\{\phi\in\KDB\,:\,\Theta\zeta_\phi=\tau\phi_0\}$, one
has that $\bar A_\Theta\phi=\A\phi_0$ for any $\phi$ in such a set and the
proof is done.
\end{proof}
\section{Singular perturbations of generalized abstract wave equations}
In this section we look for singular perturbations of operators of the kind 
$$
W_g(\phi,\psi):=(\CDB\phi+\psi,C_1\psi+\A\phi)\,.
$$
Let us begin with the simpler case in which
$$
W_g(\phi,\vp):=(\vp,C\vp+\A\phi)\,,
$$
where
$$
C:\HU\subseteq\HZ\to\HZ
$$
is a skew-adjoint operator such that:\p
H4.1)
$$
\forall\,\phi\in\HU\,,\qquad\|C\phi\|_0\le c\,\|B\phi\|_0\,;
$$
H4.2) $$C(\H_2)\subseteq\H_1\quad \text{\rm and}\quad 
\forall\phi\in\HD\,,\quad B C\phi=C B\phi\,.
$$ 
\begin{lemma} If {\rm H4.1} and {\rm H4.2} hold true then 
$$B^2-\lambda C+\lambda^2:\HD\subseteq\HZ\to\HZ$$ is invertible 
for all $\lambda\not=0$,  
$$R(\lambda):=(B^2-\lambda C+\lambda^2)^{-1}\in\BO(\HZ,\HD)\,.$$
and
$$\forall\phi\in\HU\,,\quad
\|B(B^2+\lambda^2)R(\lambda)\phi\|_0\le c\,\|B\phi\|_0\,.
$$
\end{lemma}
\begin{proof} By our hypotheses one has  
$$
\forall\,\phi\in\HD\,,\qquad\|BC\phi\|_0=\|CB\phi\|_0\le c\,\|B^2\phi\|_0\,.
$$
Thus, by induction,
$$
\forall\,k\ge 1\,,\ 
\forall\,\phi\in\H_{k+1}\,,\qquad\|B^kC\phi\|_0\le c\,\|B^{k+1}\phi\|_0\,,
$$
and $C(\H_{k+1})\subseteq\H_k$ for any $k\ge 1$. By
$$
\forall\,\phi\in\H_3\,,\qquad B^2C\phi=BCB\phi=CB^2\phi
$$
one gets 
$$\forall\,\phi\in\H_1\,,
\qquad R_0(\lambda)C\phi=CR_0(\lambda)\phi\,,
$$ 
so that $CR_0(\lambda)$ is
skew-adjoint. Thus $\uno-\lambda CR_0(\lambda) $ is boundedly invertible
for all $\lambda\not=0$ and 
$$
R(\lambda)=(\uno-\lambda CR_0(\lambda) )^{-1}R_0(\lambda)
=R_0(\lambda)(\uno-\lambda CR_0(\lambda) )^{-1}\,.
$$ 
This gives
\begin{align*}
&\|(B^2+\lambda^2)R(\lambda)\phi\|_0
=\|(\uno-\lambda CR_0(\lambda))^{-1}\phi\|_0\\
\le&\|(\uno-\lambda CR_0(\lambda))^{-1}\|_{\HZ,\HZ}\|\phi\|_0
\end{align*}
and
\begin{align*}
&\|B(B^2+\lambda^2)R(\lambda)\phi\|_0
=\|B(\uno-\lambda CR_0(\lambda))^{-1}\phi\|_0\\
=&\|(\uno-\lambda CR_0(\lambda))^{-1}B\phi\|_0
\le\|(\uno-\lambda CR_0(\lambda))^{-1}\|_{\HZ,\HZ}\|B\phi\|_0\,.
\end{align*}
\end{proof}
Let $\CB\in\BO(\HUB,\HZ)$ be the
closed bounded extension of operator
$$C:\HU\subseteq\HUB\to\HZ\,. $$
It exists by H4.1. Let $\RB\in\BO(\HUB,\HB_3)$ the closed bounded
extension of
$$
R(\lambda):\HU\subseteq\HUB\to\HB_3\,.
$$
It exists by Lemma 4.1. For such an
extension the following relations hold true:
$$
(-\A-\lambda\CB)\RB+\lambda^2\RB=\uno_{\HUB}\,,
$$
$$
R(\lambda)(-\A-\lambda\CB)+\lambda^2\RB=\uno_{\HDB}\,.
$$
Proceeding as in theorem 2.5 one obtains 
the following
\begin{theorem} Under hypotheses {\rm H4.1} and {\rm H4.2} 
the linear operator 
$$
W_g:\HDB\times\HU\subseteq\HUB\oplus\HZ\to\HUB\oplus\HZ\,,\quad
W_g(\phi,\vp):=(\vp,C\vp+\A\phi)\,,
$$
is skew-adjoint and its
resolvent is given by
\begin{align*}
&(-W_g+\lambda)^{-1}(\phi,\vp)\\
=&(\lambda\RB\phi+R(\lambda)(-\,\CB\phi+\vp),
-\phi+\lambda^2\RB\phi
+\lambda R(\lambda)(-\,\CB\phi+\vp))\,.
\end{align*}
\end{theorem}
\begin{remark} We used the notation $\HDB\times\HU$ for $D(W)$ since,
when $C\not= 0$, the scalar product inducing the graph norm on $D(W_g)$
is different from the one of $\HDB\oplus\HU$.
\end{remark}
By the previous theorem  
$$
\tilde W_g:\HDB\times\HU\times\fht\subseteq\HUB\oplus\HZ\oplus\fht
\to \HUB\oplus\HZ\oplus\fht\,,
$$
$$
\tilde W_g(\phi,\vp,\zeta):=(W_g(\psi,\vp),0)
$$
is skew-adjoint and $(-\tilde W_g+\lambda)^{-1}(\phi,\vp,\zeta)
=((-W_g+\lambda)^{-1}(\phi,\vp),\lambda^{-1}\zeta)$.\par
Now we consider a sequence $J_\nu:\HZ\to \HZ$, $\nu>0$, of self-adjoint
operators such that\smallskip\p
1. $$J_\nu\in\BO(\H_k,\H_{k+1})\,,\qquad k\ge 0\,;$$\p 
2. $$
\forall\,\phi\in\H_1\,\quad J_\nu B\phi=BJ_\nu\phi\,\quad J_\nu C\phi=CJ_\nu\phi\,;
$$
3.$$\forall\,\phi\in\HZ\,,\qquad
\lim_{\nu\downarrow
0}\, \|J_\nu\phi-\phi\|_0=0\,.
$$
Such sequence $J_\nu$ can be obtained by considering, for
example, the family 
$(\nu B^2+1)^{-1}$, but other choices are possible (see Example 3 in
the next section). We remark that the successive construction will
depend on the choice we make for such a family.\par
Denoting by $\bar J_\nu\in\BO(\HB_k,\HB_{k+1})$, $k\ge1$, the closed
bounded extension of $J_\nu$ and given $\tau\in\BO(\HDB,\fh)$ we
define 
the bounded linear map
$$
\tau_\nu:=\tau \bar J_\nu:\HUB\to\fh\,,
$$
and
$$
D(\taub):=\{\phi\in\HUB\,:\,\lim_{\nu\downarrow
0}\,\tau_\nu\phi\quad \text{\rm exists in $\fh$}\}\,,
$$
$$
\taub:D(\taub)\subseteq\HUB\to\fh\,,\quad\taub\phi:=\lim_{\nu\downarrow
0}\,\tau_\nu\phi\,.
$$
Note that for all $\phi\in\HDB$, by 3,
\begin{align*}
&\lim_{\nu\downarrow
0}\,\|B\B(\bar J_\nu\phi-\phi)\|_0^2+\|\B(\bar J_\nu\phi-\phi)\|_0^2\\
=&\lim_{\nu\downarrow
0}\,\|J_\nu B\B\phi-B\B\phi)\|_0^2+\| J_\nu\B\phi-\B\phi\|_0^2=0\,,
\,,
\end{align*}
so that $\HDB\subseteq D(\taub)$ and $\taub=\tau$ on $\HDB$.\par
Defining then
$$
\taut:D(\taub)\times\HZ\times\fht\to\fh\,,\quad\taut(\phi,\psi,\zeta):=
\taub\phi-\Theta\zeta
$$
we have that $\tau_\Theta$ satisfies A.1 and A.2.\par
Now we define $\Gb\in\BO(\HZ,\fh)$ and $\G\in\BO(\fh,\HZ)$ by
$$
\Gb:=\tau R(\lambda)\,,\qquad \G:=\breve G(-\lambda)^*\,.
$$ 
We also define $\GBb\in\BO(\HUB,\fh)$ and $\GB\in\BO(\fh,\HUB)$ by
$$
\GBb:=\tau\RB\,,\qquad\GB:=\breve{\bar G}(-\lambda)^*\,.
$$ 
Obviously $\GBb=\Gb$ on $\HU$. 
As in the previous section one has the following
\begin{lemma} $$\ran(\G)\subseteq\ran(B)\qquad\text{\it and}\qquad 
\GB=\B^{-1}B^{-1}\G\,.$$
$$\ran(\GB)\subseteq\HDB\qquad\text{\it and}\qquad\GB\in\BO(\fh,\HDB)\,.$$
\end{lemma}
Now we define
$\GTb\in\BO(\HUB\oplus\HZ\oplus\fht,\fh)$ by
\begin{align}
&\GTb(\phi,\vp,\zeta):=
\taut(-\tilde W+\lambda)^{-1}(\phi,\vp,\zeta)\\
=&\lambda\GBb\phi+\Gb(-\,\CB\phi+\vp)-\lambda^{-1}\Theta \zeta\,
\end{align}
and
$\GT\in\BO(\fh,\HUB\oplus\HZ\oplus\fht)$ by
\begin{align}
&\GT \zeta:=-\,\breve G_\Theta(-\lambda)^*\zeta\\
=&(\lambda\GB \zeta+\CB^*\G\zeta,-\,\G \zeta,-\lambda^{-1}\zeta)\,.
\end{align}
Regarding the adjoint of $\bar C$ one has the following
\begin{lemma} 
$$
\CB^*=-\B^{-1}\CB\B^{-1}\,.
$$
\end{lemma}
\begin{proof} Since $C$ commutes with $B$, for any $\phi$ in $\HU$ one has 
$$
\CB\B^{-1}\phi=\B^{-1}C\phi\,,
$$
and thus for any $\phi$ and $\psi$
in $\HU$ one has
\begin{align*}
&[-\B^{-1}\CB\B^{-1}\phi,\psi]_1=-\langle\CB\B^{-1}\phi,B\psi\rangle_0\\=
&-\langle\B^{-1}C\phi,B\psi\rangle_0=-\langle C\phi,\psi\rangle_0=
\langle \phi,C\psi\rangle_0\,.
\end{align*}
\end{proof}
Let us now consider the bounded linear map
$$
\Gamma_{\nu,\Theta}(\lambda):=
-\tau_\nu\B^{-1}(-C+\lambda) B^{-1}\G-\frac{1}{\lambda}\,\Theta\,.
$$
We have the following
\begin{lemma}
$$\Gamma_{\nu,\Theta}(\lambda)^*=-\Gamma_{\nu,\Theta}(-\lambda)\,.$$
\end{lemma}
\begin{proof} At first let us observe that, being $CR_0(\epsilon)$
skew-adjoint (see the proof of Lemma 4.1), one has
$$
\forall\,\epsilon>0\,,\qquad \|(1\pm\epsilon CR_0(\epsilon))^{-1}\|_{\HZ,\HZ}\le
1\,.
$$ 
Thus, using H4.1, functional calculus and dominated convergence theorem,
\begin{align*}
&\lim_{\epsilon\downarrow 0}\,\|((\uno\pm\epsilon
CR_0(\epsilon))^{-1}-\uno)\phi\|^2\le\lim_{\epsilon\downarrow 0}\,
\|\epsilon CR_0(\epsilon)\phi\|^2\\
\le& c^2\lim_{\epsilon\downarrow 0}\, 
\|\epsilon BR_0(\epsilon)\phi\|^2
=\lim_{\epsilon\downarrow 0}\,
\int_\RE d\mu_\phi(x)\, \left|\frac{\epsilon x}{x^2+\epsilon^2}\right|^2=0\,.
\end{align*}
Since $$
B^2R(\pm\epsilon)=B^2R_0(\epsilon)(1\pm\epsilon CR_0(\epsilon))^{-1}\,,$$
and $B^2R_0(\epsilon)$ strongly converges to $\uno_{\HZ}$ when
$\epsilon\downarrow 0$ (see the proof of Lemma 3.2), one has that 
$$
\text{\rm s-}\lim_{\epsilon\downarrow 0}B^2R(\pm\epsilon)=\uno_{\HZ}\,.
$$
This implies (proceeding as in the proof of Lemma 4.1) that
$R(\epsilon)\G$ strongly converges in $\BO(\fh,\HDB)$ to $\B^{-1}
B^{-1}\G$ when
$\epsilon\downarrow 0$. Therefore
\begin{align*}
&\Gamma_{\nu,\Theta}(\lambda)\\
&=\text{\rm s-}\lim_{\epsilon\downarrow 0}\,
-\frac{1}{\lambda}\,\left(\Theta+
\lambda\,\tau_\nu(-C+\lambda\pm\epsilon)\,(B^2\mp\epsilon
C+\epsilon^2)^{-1}\,G(\lambda)\right)\,,
\end{align*}
and the proof is concluded by showing that 
\begin{align*}
&(\tau_\nu\, C\,(B^2-\epsilon
C+\epsilon^2)^{-1}\,G(\lambda))^*\\
=&-\tau_\nu\, C\,(B^2+\epsilon
C+\epsilon^2)^{-1}\,G(-\lambda)\,.
\end{align*}
Proceeding as in \cite{[P1]}, Lemma 2.1, by 
first resolvent identity one obtains 
\begin{align*}
&(-C+\lambda+\epsilon)\,(B^2-\epsilon
C+\epsilon^2)^{-1}\,G(\lambda)\\
=&\frac{G(\epsilon)-G(\lambda)}{\lambda-\epsilon}\\
=&(-C+\lambda+\epsilon)\,(B^2-\lambda
C+\lambda^2)^{-1}\,G(\epsilon)
\,,
\end{align*}
so that
$$
(B^2-\epsilon
C+\epsilon^2)^{-1}\,G(\lambda)
=(B^2-\lambda
C+\lambda^2)^{-1}\,G(\epsilon)
\,.
$$
Therefore we need to show that
\begin{align*}
&(\tau_\nu\, C\,(B^2-\lambda
C+\lambda^2)^{-1}\,G(\epsilon))^*\\
=&-\tau_\nu\, C\,(B^2+\epsilon
C+\epsilon^2)^{-1}\,G(-\lambda)\,.
\end{align*}
Since $C$, $B$ and $J_\nu$ commute, we have
\begin{align*}
&(\tau_\nu\, C\,(B^2-\lambda
C+\lambda^2)^{-1}\,G(\epsilon))^*
=(\breve G(\lambda)\, C\,J_\nu\,\,G(\epsilon))^*\\
=&G(\epsilon)^*(C\,J_\nu)^*\breve G(\lambda)^*
=-\breve G(-\epsilon)\,C\,J_\nu\,G(-\lambda)\\
=&-\tau_\nu\, C\,(B^2+\epsilon
C+\epsilon^2)^{-1}\,G(-\lambda)
\end{align*}
and the proof is done.
\end{proof}

Now we suppose that \p
H4.3)
$$\ran(\CB^*\G)\equiv\ran(\CB\,\GB)\subseteq
D(\taub)\,.$$
Note that H4.3 is always verified if $\CB\in\BO(\HDB,\HDB)$, but such a hypothesis can 
hold true also in situations where $C$ is unbounded (see Example 3 in the next section).
Then, by the uniform boundedness principle 
$\taub\,\CB^*\G\in\BO(\fh)$, 
so that A.4 is satisfied,
\begin{align}
\GMT:=&-\,\taut\GT
=-\lambda\tau\B^{-1} B^{-1}\G-\taub\CB^*\G-\frac{1}{\lambda}\,\Theta\\
=&-\taub\B^{-1}(-C+\lambda) B^{-1}\G-\frac{1}{\lambda}\,\Theta
\end{align}
is well-defined and $\GMT\in\BO(\fh)$. By H4.3 one has 
$$
\GMT=\text{\rm s-}\lim_{\epsilon\downarrow 0}\Gamma_{\nu,\Theta}(\lambda)\,,
$$ 
so that the previous lemma implies
$$\Gamma_\Theta(\lambda)^*=-\Gamma_\Theta(-\lambda)\,.$$
Thus A.5 is satisfied. Suppose now that H3.1 holds true. 
Then, since $$R(\lambda)=
R_0(\lambda)(\uno-\lambda CR_0(\lambda) )^{-1}\,,$$ 
one has
$$
\ran(\G)\cap\HU=\left\{0\right\}\,,
$$
so that A.3 is satisfied. In conclusion, by Theorem 6.2
we can define a skew-adjoint extension of the skew-symmetric operator 
given by restricting $\tilde W$ to the dense set
$$\N_\Theta:=\left\{(\phi,\vp,\zeta)\in \HDB\times\HU\times\fht\ :\
\tau\phi=\Theta \zeta\right\}\,:
$$
\begin{theorem} Suppose that \text{\rm H3.0, H3.1, H4.1, H4.2} and \text{\rm
H4.3} hold true. Let 
\begin{align*}
&D(\tilde W_\Theta)\\:=&\left\{(\phi_0,\vp_0,\zeta_\phi)\in
\HUB\oplus\HZ\oplus\fht\ :\ \phi_0=\phi_\lambda+\B^{-1}CB^{-1}\G
\zeta_\vp\,,\right.\\
&\left.\vp_0=\vp_\lambda+\G \zeta_\vp\,,
\ \phi_\lambda\in\HDB\,,\ \vp_\lambda\in\HU\,,\  
\zeta_\vp\in\fh\,,\ \Theta \zeta_\phi=\taub\phi_0\,\right\}\,.
\end{align*}
Then
$$
\tilde W_\Theta:D(\tilde W_\Theta)\subseteq 
\HUB\oplus\HZ\oplus\fht\to
\HUB\oplus\HZ\oplus\fht\,,$$
$$
\tilde W_\Theta(\phi_0,\vp_0,\zeta_\phi)
:=(\vp_\lambda-\lambda\B^{-1}(-C+\lambda)B^{-1}\G
\zeta_\vp,\, C\vp_\lambda+\A\phi_\lambda,\zeta_\vp)
$$
is a skew-adjoint extension of the restriction of 
$$\tilde W_g:\HDB\times\HU\times\fht\subseteq\HUB\oplus\HZ\oplus\fht\to
\HUB\oplus\HZ\oplus\fht\,,$$
$$
\tilde W_g(\phi,\psi,\zeta)=(\psi,C\phi+\A\phi,0)$$ 
to the dense set $\N_\Theta$. The resolvent of $\tilde W_\Theta$ is given by
$$
(-\tilde W_\Theta+\lambda)^{-1}=(-\tilde
W_g+\lambda)^{-1}+\GT\GMT^{-1}\GTb\,,
$$
where the linear operators $(-\tilde
W+\lambda)^{-1}$, $\GTb$, $\GT$, $\GMT$,  have been defined in
\text{\rm Theorem
4.2, (4.2), (4.4)} and \text{\rm (4.6)} respectively. 
\end{theorem}
\begin{proof} By Theorem 6.4 we known that 
$(-\tilde W_g+\lambda)^{-1}+\GT\GMT^{-1}\GTb$ is the resolvent
of a skew-adjoint extension $\hat W_\Theta$ of the restriction of 
$\tilde W_g$ to 
the dense set $\N_\Theta$. Therefore $(\hat \phi_0,\hat
\vp,\hat \zeta_\phi)\in D(\hat W_\Theta)$
if and only if
\begin{align*}
\hat\phi_0=&\,\hat\phi_\lambda+
\left(\lambda\B^{-1}B^{-1}+\CB^*\right)\G
\GMT^{-1}(\tau\hat\phi_\lambda-\Theta
\zeta_\lambda)\,,\\
\hat\vp=&\,\vp_\lambda-\G\GMT^{-1}(\tau\hat\phi_\lambda-\Theta
\zeta_\lambda)\,,\\
\hat\zeta_\phi=&\,\zeta_\lambda
-\frac{1}{\lambda}\,\GMT^{-1}(\tau\hat\phi_\lambda-\Theta
\zeta_\lambda)\,,
\end{align*} 
where
$$
\hat\phi_\lambda\in\HDB\,,\qquad \vp_\lambda\in\HU\,,
\qquad \zeta_\lambda\in\fh\,.
$$
Let us now show that $D(\hat W_\Theta)=D(\tilde W_\Theta)$. \par
Since $\ran(\GB)\subseteq\HDB$, so that 
$$
\hat\phi_\lambda+\lambda\B^{-1}B^{-1}\G\GMT^{-1}(\tau\hat\phi_\lambda-\Theta
\zeta_\lambda)\in\HDB\,,
$$
and
\begin{align*}
&\taut((-\hat W_g+\lambda)^{-1}+\GT\GMT^{-1}\GTb)\\
&=\GTb-\GMT\GMT^{-1}\GTb=0\,, 
\end{align*}
so that $\tau\hat\phi_0=\Theta \hat \zeta_\phi$, we have $D(\hat
W_\Theta)\subseteq D(\tilde W_\Theta)$. Let us now prove the reverse
inclusion. Given $(\phi_0,\vp_0,\zeta_\phi)\in D(\tilde W_\Theta)$ let
us define 
$$
\hat\phi_\lambda:=\,\phi_\lambda+\lambda\GB\zeta_\vp\,,\qquad
\zeta_\lambda:=\,\zeta_\phi-\frac{1}{\lambda}\zeta_\vp\,.
$$
Then
$$
\taub\phi_0
=\tau\hat\phi_\lambda-\lambda\tau\GB\zeta_\vp
-\taub\CB^*\G \zeta_\vp
=\Theta \zeta_\phi=\Theta\left(\zeta_\lambda
+\frac{1}{\lambda}\,\zeta_\vp\right)
$$
implies
$$
\tau\hat\phi_\lambda-\Theta \zeta_\lambda=\left(
\lambda\tau\GB+\taub\CB^*\G+\frac{1}{\lambda}\,\Theta\right)\zeta_\vp\,,
$$
i.e.
$$
\zeta_\vp
=-\GMT^{-1}(\tau\hat\phi_\lambda-\Theta \zeta_\lambda)\,.
$$
Thus $D(\tilde W_\Theta)\subseteq D(\hat W_\Theta)$. Now we have
\begin{align*}
&\hat W_\Theta(\phi_0,\psi,\zeta_\phi)=\tilde
W_g(\hat\phi_\lambda,\vp_\lambda,\zeta_\lambda)+\lambda
(\phi_0-\hat\phi_\lambda,\vp-\vp_\lambda,\zeta_\phi-\zeta_\lambda)\\
=&(\vp_\lambda
-\lambda^2\GB\zeta_\vp-\lambda\CB^*\G \zeta_\vp,
C\vp_\lambda+\A\hat\phi_\lambda+\lambda\G\zeta_\vp,\,\zeta_\vp)\\
=&(\vp_\lambda
-\lambda^2\GB\zeta_\vp-\lambda\CB^*\G \zeta_\vp,
C\vp_\lambda+\A(\phi_\lambda+\lambda\GB\zeta_\vp)
+\lambda\G \zeta_\vp,\,\zeta_\vp)\\
=&(\vp_\lambda-\lambda(\lambda\B^{-1}B^{-1}-\B^{-1}CB^{-1})\G \zeta_\vp,
C\vp_\lambda+\A\phi_\lambda,\,\zeta_\vp)\\
=&(\vp_\lambda-\lambda\B^{-1}(-C+\lambda)B^{-1}\G
\zeta_\vp,\, C\vp_\lambda+\A\phi_\lambda,\zeta_\vp)\\
=&\tilde W_\Theta(\phi_0,\psi,\zeta_\phi)\,.
\end{align*}
\end{proof}
Let us now consider two skew-adjoint operators
$$
C_1:\HU\subseteq\HZ\to\HZ\,,\qquad C_2:\HU\subseteq\HZ\to\HZ\,,
$$
such that \p
H4.1.1)
$$
\forall\,\phi\in\HU\,,\quad\|C_1\phi\|_0\le c_1\,\|B\phi\|_0\,,\quad
\|C_2\phi\|_0\le c_2\,\|B\phi\|_0\,,\quad c_1c_2<1
$$
H4.2.1)
$$C_1(\HD)\subseteq\HU\,,\qquad C_2(\HD)\subseteq\HU\,.$$ 
and 
$$
\forall\phi\in\HD\,,\quad C_1C_2\phi=C_2C_1\phi\,,\quad
BC_1\phi=C_1B\phi\,,\quad
BC_2\phi=C_2B\phi\,.
$$
Then by the Kato-Rellich theorem
$$
-A_C:=B^2+C_1C_2:\HD\subseteq\HZ\to\HZ
$$
is self-adjoint, positive and injective. Let $B_C$ be the
self-adjoint, positive and injective operator defined by 
$B_C:=(-A_C)^{1/2}$. Since, by H4.2.1,  
$$
(1-c_1c_2)\|B\phi\|\le \|B_C\phi\|\le (1+c_1c_2)\|B\phi\|\,,
$$ 
the domain of $B_C$ coincides with the space $\HU$, the domain of
$B$. 
Moreover, since $B$ and $B_C$ commutes, 
$$
(1-c_1c_2)^k\|B^k\phi\|\le \|B^k_C\phi\|\le (1+c_1c_2)^k\|B^k\phi\|\,,
$$ 
thus the Hilbert spaces generated by $B_C$ coincide, 
as Banach spaces (in the sense that each space has an equivalent
norm), with the ones
generated by $B$, i.e. coincide with $\H_k$, $\bar\H_k$, and
$\H_{-k}$, $k\ge 1$.\par
Let $\A_C:=-B_C\B_C\in\BO(\HDB,\HZ)$, where $\B_C\in\BO(\HUB,\HZ)$ 
is the closed bounded extension of $B_C:\HU\subseteq\HUB\to\HZ$. We
know that $\A_C$ coincides with the closed bounded extension of
$A_C:\HD\subseteq\HDB\to\HZ$. Since $C_2$ commutes with $B$,
by H4.2.1 we have 
$$
\|BC_2\phi\|_0\le c_2\|B^2\phi\|_0\,.
$$ 
Thus we can define $\CDB\in\BO(\HUB,\HZ)\cap\BO(\HDB,\HU)$ 
as the closed bounded
extension of $C_2:\HU\subseteq\HUB\to\HZ$ and 
$$
\A_C=\A-C_1\CDB\,.
$$
Since $C:=C_1+C_2$ and $B_C:=\sqrt{B^2+C_1C_2}$ satisfy H4.1 and H4.2,
by Theorem 4.2 we have that
$$W_g:\HDB\times\HU\subseteq\HUB\oplus\HZ\to\HUB\oplus\HZ$$
$$
W_g(\phi,\vp):=(\vp,(C_1+C_2)\vp+(\A-C_1\CDB)\phi)
$$
is skew-adjoint once we put on $\HUB\oplus\HZ$ the scalar product 
$$
\langle\langle\,(\phi_1,\vp_1),(\phi_2,\vp_2)\,\rangle\rangle
:=\langle\bar B_C\phi_1,\bar B_C\phi_2\rangle_0
+\langle\vp_1,\vp_2\rangle_0\,.
$$
Let us define the Hilbert space 
$(\H_C,\langle\langle\cdot,\cdot\rangle\rangle_C)$ by
$\H_C=\HUB\times\HZ$, 
\begin{align*}
&\langle\langle\,(\phi_1,\psi_1),(\phi_2,\psi_2)\,\rangle\rangle_C\\
:=&\langle\bar B_C\phi_1,\bar B_C\phi_2\rangle_0
+\langle\psi_1+\CDB\phi_1,\psi_2+\CDB\phi_2\rangle_0\\
=&\langle\bar B\phi_1,\bar B\phi_2\rangle_0
+\langle\CDB\phi_1,\psi_2\rangle_0+
\langle\psi_1,\CDB\phi_2\rangle_0+\langle\psi_1,\psi_2\rangle_0\\
&+\langle(\CDB-\CUB)\phi_1,\CDB\phi_2\rangle_0\,,
\end{align*}
where $\CUB\in\BO(\HUB,\HZ)$ denotes the closed bounded
extension of $$C_1:\HU\subseteq\HUB\to\HZ\,.$$ 
Then the map 
$$
S:\HUB\oplus\HZ\to\H_C\,,\quad S(\phi,\vp):=(\phi, \vp-\CDB\phi)
$$
is unitary
and the linear operator
$$
SW_gS^*:\HDB\times\HU\subseteq\H_C\to\H_C
$$
\begin{align*}
&SW_gS^*(\phi,\psi)=SW_g(\phi,\psi+\CDB\phi)\\
=&
S(\CDB\phi+\psi,(C_1+C_2)(\psi+\CDB\phi)+(\A-C_1\CDB)\phi)\\
=&
S(\CDB\phi+\psi,(C_1+C_2)\psi+(\A+C_2\CDB)\phi)\\
=&
(\CDB\phi+\psi,(C_1+C_2)\psi+(\A+C_2\CDB)\phi-
\CDB(\CDB\phi+\psi))\\
=&(\CDB\phi+\psi,C_1\psi+\A\phi)
\end{align*}
is skew-adjoint.\par
Let us now define, on the Hilbert space $\H_C\oplus\fht$ with
scalar product 
$$
\langle\langle\,(\phi_1,\psi_1,\zeta_1),
(\phi_1,\psi_2,\zeta_2)\,\rangle\rangle_{C,\Theta}:=
\langle\langle\,(\phi_1,\psi_1),
(\phi_1,\psi_2)\,\rangle\rangle_C+\langle\zeta_1,\zeta_2\rangle_\Theta\,,
$$
the skew-adjoint operator
$$
\tilde W_g:\HDB\times\HU\times\fht\subseteq\H_C\oplus\fht
\to \H_C\oplus\fht\,,
$$
$$
\tilde W_g(\phi,\psi,\zeta):=(\CDB\phi+\psi,C_1\psi+\A\phi,0)\,.
$$
Let $$\GCb:=\tau
(B^2+(-C_1+\lambda)(-C_2+\lambda))^{-1}\,,\quad \GC:=\breve
G_C(-\lambda)^*\,,$$ 
and suppose\p
H4.3.1)
$$\ran(\CB_j^*G_C(\lambda))\subseteq
D(\taub)
$$
where now 
$$
\CB_j^*=-\bar B_C^{-1}\CB_j\bar B_C^{-1}\,,\quad j=1,2\,,
$$
Then by the previous theorem we obtain the following
\begin{theorem} Suppose \text{\rm H3.0, H3.1, H4.1.1, H4.2.1} and \text{\rm
H4.3.1} hold true. Then the linear operator 
$$
\tilde W_\Theta:D(\tilde W_\Theta)\subseteq 
\H_C\oplus\fht\to
\H_C\oplus\fht\,,$$
\begin{align*}
&D(\tilde W_\Theta):=\{(\phi_0,\psi_0,\zeta_\phi)\in
\HUB\times\HZ\times\fht\, :\,\\
& \phi_0=\phi_\lambda+\B_C^{-1}(C_1+C_2)B_C^{-1}\GC
\zeta_\psi\,,\\
&\psi_0=\psi_\lambda+(\uno-\CDB\B_C^{-1}(C_1+C_2)B_C^{-1})\GC \zeta_\psi\,,\\
&\phi_\lambda\in\HDB\,,\ \psi_\lambda\in\HU\,,\  
\zeta_\psi\in\fh\,,\ \Theta \zeta_\phi=\taub\phi_0\,\}\,,
\end{align*}
\begin{align*}
&\tilde W_\Theta(\phi_0,\psi_0,\zeta_\phi)\\
:=(&\CDB\phi_\lambda
+\psi_\lambda-\lambda\B_C^{-1}(-(C_1+C_2)+\lambda)B_C^{-1}\GC
\zeta_\psi,\\ 
&C_1\psi_\lambda+\A\phi_\lambda
+\lambda\CDB\B_C^{-1}(-(C_1+C_2)+\lambda)B_C^{-1}\GC\zeta_\psi,\zeta_\psi)
\end{align*}
is a skew-adjoint extension of the restriction of $$
\tilde W_g:\HDB\times\HU\times\fht\subseteq\H_C\oplus\fht
\to \H_C\oplus\fht\,,
$$
$$
\tilde W_g(\phi,\psi,\zeta):=(\CDB\phi+\psi,C_1\psi+\A\phi,0)\,.
$$
 to 
the dense set $\N_\Theta$. 
\end{theorem}
Let us now suppose that H3.2 holds true. Then we can define 
$$
G_C:\fh\to\V\,,\qquad G_C
:=\GC+\lambda\B_C^{-1}(-(C_1+C_2)+\lambda)B_C^{-1}\GC\,.
$$
\begin{lemma} The definition of $G_C$ is $\lambda$-independent. Moreover 
$$
\ran(G_C)\cap\HUB=\left\{0\right\}\,.
$$
\end{lemma}
\begin{proof} Let $C=C_1+C_2$. Proceeding as in \cite{[P1]}, Lemma 2.1, by 
first resolvent identity one obtains 
$$
(\lambda-\mu)\,(-C+\lambda+\mu)\,(B_C^2-\mu
C+\mu^2)^{-1}\,G_C(\lambda)=
G_C(\mu)-G_C(\lambda)\,,
$$
i.e.
\begin{align*}
&(\uno+\mu\bar B_C^{-1}(-C+\mu)B_C^{-1})\,
((\lambda-\mu)\,(-C+\lambda+\mu))\times\\
&\times(B_C^2-\mu
C+\mu^2)^{-1}\,G_C(\lambda)\\
=&\bar B_C^{-1}
((\lambda-\mu)\,(-C+\lambda+\mu))B_C^{-1}G_C(\lambda)\\
=&(\uno+\mu\bar B_C^{-1}(-C+\mu)B_C^{-1})\,
(G_C(\mu)-G_C(\lambda))\,.
\end{align*}
This implies 
\begin{align*}
&\GC+\lambda\B_C^{-1}(-C+\lambda)B_C^{-1}\GC\\
=&G_C(\mu)+\mu\B_C^{-1}(-C+\mu)B_C^{-1}G_C(\mu)\,.
\end{align*}
Suppose there exists $\zeta\in\fh$ such that
$$
\GC\zeta+\lambda\B_C^{-1}(-C+\lambda)B_C^{-1}\GC\zeta=\phi\in\HUB\,.
$$
Then $\GC\zeta\in\HU$ and so, by H3.1, $\GC\zeta=0$. Thus 
$$\B_C^{-1}(-C+\lambda)B_C^{-1}\GC\zeta=0$$ and the proof is done.
\end{proof}

For any $k\ge 0$, $j=1,2$, let 
$$\hat B:\H_{-k}\to\H_{-(k+1)}\,,$$
$$\hat B_C:\H_{-k}\to\H_{-(k+1)}\,,$$
$$ \hat C_j:\H_{-k}\to\H_{-(k+1)}\,,$$
be the closed bounded extensions of
$$B:\HU\subseteq\H_{-k}\to\H_{-(k+1)}\,,$$ 
$$B_C:\HU\subseteq\H_{-k}\to\H_{-(k+1)}$$ and 
$$\hat C_j:\HU\subseteq\H_{-k}\to\H_{-(k+1)}\,,$$ respectively. 
Define also 
$$\hat A:\HUB\to\H_{-1}\,,\qquad\hat A:=-\hat B\B$$
and 
$$\hat C_2G_C:\fh\to \H_{-1}\,,$$
$$
\hat C_2G_C:=\hat C_2\GC+\lambda\CD\B_C^{-1}(-(C_1+C_2)+\lambda)B_C^{-1}\GC\,.
$$
Then 
\begin{lemma}
$$
\tilde W_\Theta(\phi_0,\psi_0,\zeta_\phi)=
(\hat C_2\phi_0+\psi_0-G_C\zeta_\psi, \hat C_1\psi_0+\hat A\phi_0+\hat C_2G_C\zeta_\psi,\zeta_\psi)\,.
$$ 
\end{lemma}
\begin{proof} Since
$$
\hat B_CC_1=\hat C_1B_C\,,\qquad\hat B_CC_2=\hat C_2B_C\,,
$$
and
$$
-\hat A+ \hat C_1\CDB=\hat B_C\bar B_C\,, 
$$
one has
\begin{align*}
&\CDB\phi_\lambda+\psi_\lambda
-\lambda\B_C^{-1}(-(C_1+C_2)+\lambda)B_C^{-1}\GC\zeta_\psi\\
=&\hat C_2\phi_0-\CDB\B_C^{-1}(-(C_1+C_2)+\lambda)B_C^{-1}\GC\zeta_\psi\\
&+\psi_0-(\uno+\CDB\B_C^{-1}(-(C_1+C_2)+\lambda)B_C^{-1})\GC\zeta_\psi\\
&-\lambda\B_C^{-1}(-(C_1+C_2)+\lambda)B_C^{-1}\GC\zeta_\psi\\
=&\hat C_2\phi_0+\psi_0-G_C\zeta_\psi
\end{align*}
and
\begin{align*}
&C_1\psi_\lambda+\A\phi_\lambda
+\lambda\CDB\B_C^{-1}(-(C_1+C_2)+\lambda)B_C^{-1}\GC\zeta_\psi\\
=&\hat C_1\psi_0-\hat C_1(\uno-\CDB\B_C^{-1}(C_1+C_2)B_C^{-1})\GC
\zeta_\psi\\
&+\hat A\phi_0-\hat A\B_C^{-1}(C_1+C_2)B_C^{-1}\GC\zeta_\psi\\
&+\hat C_2G_C\zeta_\psi-\hat C_2\GC\zeta_\psi\\
=&\hat C_1\psi_0-\hat C_1\GC\zeta_\psi-\hat C_2\GC\zeta_\psi
+\hat A\phi_0+\hat C_2G_C\zeta_\psi\\
&(-\hat A+\hat C_1\CDB)\B_C^{-1}(C_1+C_2)B_C^{-1}\GC
\zeta_\psi\\
=&\hat C_1\psi_0+\hat A\phi_0+\hat C_2G_C\zeta_\psi
\,.
\end{align*}
\end{proof}
By Lemma 4.8 we can define the Hilbert space $(\H_\Theta, 
\langle\langle\cdot,\cdot\rangle\rangle_{\H_\Theta})$ by
\begin{align*}
&\H_\Theta:=\{(\phi,\psi)\in\V\times\H_{-1}\,:\, 
\phi=\phi_0+G_C\zeta_\phi,\\ &\psi=\psi_0-\hat C_2G_C\zeta_\phi, \ 
(\phi_0,\psi_0,\zeta_\phi)\in\H_C\oplus\fh\}
\end{align*}
with scalar product
$$
\langle\langle\,(\phi,\psi),(\tilde\phi,\tilde\psi)\,\rangle\rangle_{\H_\Theta}:=
\langle\langle\,(\phi_0,\psi_0,\zeta_{\phi}),
(\tilde\phi_0,\tilde\psi_0,\zeta_{\tilde\phi})\,\rangle\rangle_{C,\Theta}\,,
$$ 
so that map 
$$
U:\H_C\oplus\fh_\Theta\to\H_\Theta\,,\qquad 
U(\phi_0,\psi_0,\zeta_\phi):=(\phi_0+G_C\zeta_\phi,\psi_0-\hat C_2G_C\zeta_\phi)
$$
is unitary. Thus in conclusion we have the following
\begin{theorem} Suppose \text{\rm H3.0, H3.1, H3.2, H4.1.1, H4.2.1} and \text{\rm
H4.3.1} hold true. Then the linear operator 
$$
W_\Theta:D(W_\Theta)\subseteq\H_\Theta\to\H_\Theta\,,
$$
\begin{align*}
D(W_\Theta):=\{&(\phi,\psi)\in\H_\Theta\, :\,
 \phi_0=\phi_\lambda+\B_C^{-1}(C_1+C_2)B_C^{-1}\GC
\zeta_\psi\,,\\
&\psi_0=\psi_\lambda+(\uno-\CDB\B_C^{-1}(C_1+C_2)B_C^{-1})\GC \zeta_\psi\,,\\
&\phi_\lambda\in\HDB\,,\ \psi_\lambda\in\HU\,,\  
\zeta_\psi\in\fh\,,\ \Theta \zeta_\phi=\taub\phi_0\,\}\,.
\end{align*}
$$
W_\Theta(\phi,\psi):=U \tilde W_\Theta U^*(\phi,\psi)=
(\hat C_2\phi_0+\psi_0,\hat C_1\psi_0+\hat A\phi_0)\,,
$$
is skew-adjoint. It coincides with 
$$W_g:\HDB\times\HU\subseteq\H_C\to\H_C$$
$$
W_g(\phi,\psi):=(\CDB\phi+\psi,C_1\psi+\A\phi)
$$
on the dense set
$$D(W_g)\cap D(W_\Theta)=\left\{\phi\in\HDB\ :\
\tau\phi=0\right\}\times\HU\,. $$
\end{theorem}

\section{examples}

\p
{\it -- Example 1.} Let $A_0$ be the negative and injective self-adjoint
operator on ${\mathcal H}_0=L^2(0,\infty)$ corresponding to the second
derivative operator with Dirichlet boundary contidions at zero, i.e. 
$$
A_0:D(A_0)\subset L^2(0,\infty)\to L^2(0,\infty)\,,\quad A_0\phi:=\phi''\,,
$$
where $\H_2\equiv D(A_0)\equiv H^2_0(0,\infty)$,
$$
H^2_0(0,\infty)
:=\left\{\phi\in L^2(0,\infty)\,:\,\phi''\in L^2(0,\infty)\,,\ 
\phi(0_+)=0\right\}.
$$
Let $B_0$ be the positive and injective self-adjoint
operator defined by $B_0:=\sqrt{-A_0}$. We have 
$\H_1\equiv D(B_0)\equiv H^1_0(0,\infty)$, where
$$
H^1_0(0,\infty):=\left\{\phi\in L^2(0,\infty)\,:\,\phi'\in L^2(0,\infty)\,,\ \phi(0_+)=0\right\}\,,
$$
with scalar product 
$$\langle\phi_1,\phi_2\rangle_{1}:=\langle\phi_1,\phi_2\rangle+
\langle\phi_1',\phi_2'\rangle\,.
$$  
Here $\langle\cdot,\cdot\rangle$ denotes here the usual scalar
product on $L^2(0,\infty)$.\par
Let us now consider $\bar\H_1\equiv\bar H^1_0(0,\infty)$, the completion of 
$H^1_0(0,\infty)$ with respect to the scalar product 
$$[\phi_1,\phi_2]_{1}:=
\langle\phi_1',\phi_2'\rangle\,.
$$   
One has 
$$
\bar H^1_0(0,\infty):=\left\{\phi\in \bigcup_{b>0} L^2(0,b)\,:\,\phi'\in L^2(0,\infty)\,,\ \phi(0_+)=0\right\}\,,
$$
and then 
$$
\bar H^2_0(0,\infty)
:=\left\{\phi\in \bigcup_{b>0} L^2(0,b)\,:\,\phi'\,,\phi''\in
L^2(0,\infty)\,,\ \phi(0_+)=0\right\}.
$$
Moreover $\bar A_0$ acts on $\bar H^2_0(0,\infty)$ as the
second (distributional) derivative operator. The resolvent 
$(-A_0+\lambda^2)^{-1}$ has an integral kernel given by 
$$
\g_D(\lambda;x,y)=\frac{e^{-|\lambda|\,|x-y|}-e^{-|\lambda|\,(x+y)}}{ 2|\lambda|}\,.
$$
We consider now the negative and injective self-adjoint operator on
$\bigoplus_{k=1}^n L^2(0,\infty)$ defined by
$A:=\bigoplus_{k=1}^n A_0$ and the bounded linear map 
$$
\tau:\bigoplus_{k=1}^n\bar H^2_0(0,\infty)\to\C^n\,,\quad 
\tau(\phi_1,\dots,\phi_n):=\left(\phi_1'(0_+),\dots,\phi_n'(0_+)\right)\,.
$$
Obviously
$\tau_\Theta(\phi_1,\dots,\phi_n,\zeta):=\tau(\phi_1,\dots,\phi_n)
-\Theta\zeta$ satisfies hypothesis H.3.0 for any positive and
injective Hermitean $\Theta$.\par
One has that $G(\lambda):\C^n\to \bigoplus_{k=1}^n L^2(0,\infty)$ is
represented by the vector in $\bigoplus_{k=1}^n L^2(0,\infty)$ given by
$$
\g_\lambda(x_1,\dots,x_n)
:=\left(e^{-|\lambda|{x_1}},\dots,e^{-|\lambda|{x_n}}\right)\,,
$$
while $\bar G(\lambda):\C^n\to \bigoplus_{k=1}^n \bar H^2_0(0,\infty)$ is
represented by the vector in $\bigoplus_{k=1}^n \bar H^2_0(0,\infty)$ given by
\begin{align*}
&\bar\g_\lambda(x_1,\dots,x_n)\\
\equiv&\lim_{\epsilon\to 0}\, \left(
\int_0^\infty dx_1\,\g_D(\epsilon;x_1,y_1)\,\g_\lambda(y_1),\dots,
\int_0^\infty dx_n\,\g_D(\epsilon;x_n,y_n)\,\g_\lambda(y_n)
\right)\\
=&\left(\frac{e^{-|\lambda|{x_1}}-1}{|\lambda|^2},\dots,
\frac{1-e^{-|\lambda|{x_1}}}{|\lambda|^2}\right)\,.
\end{align*} 
Therefore
$$
\Gamma_\Theta(\lambda)=-\lambda\tau\bar
G(\lambda)-\frac{1}{\lambda}\,\Theta=-\frac{1}{\lambda}\,(|\lambda|+\Theta)\,.
$$
Note that, since $\g_\lambda(0_+)\not=0$, Ran$(G(\lambda))\cap
H^1_0(0,\infty)=\left\{0\right\}$ and H.3.1 is satisfied.\par
Hypothesis H.3.2 is satisfied by taking  ${\mathcal V}
=\bigcup_{b>0}L^2((0,b)^n)$ and
$G:\C\to\bigcup_{b>0}L^2((0,b)^n)$ is represented by the
constant vector 
$$\g(x_1,\dots,x_n)=\left(\g_\lambda+\lambda^2\bar\g_\lambda\right)(x_1,\dots,x_n)
=(1,\dots,1)\,.$$ 
Defining
$$
\bar H^1(0,\infty):=\left\{\phi\in \bigcup_{b>0} L^2(0,b)\ :\ 
\phi'\in L^2(0,\infty)\right\}\,,
$$
$$
H^1(0,\infty):=\bar H^1(0,\infty)\cap L^2(0,\infty)\,,
$$ 
and
$$
\bar H^2(0,\infty):=\left\{\phi\in \bigcup_{b>0} L^2(0,b)\ :\ 
\phi'\,,\phi''\in L^2(0,\infty)\right\}\,,
$$
$$
H^2(0,\infty):=\bar H^2(0,\infty)\cap L^2(0,\infty)\,,
$$
one has 
\begin{align*}
\bar K^1:=&
\left\{\Phi=\Phi_0+\zeta_\Phi\g\,,\ \Phi_0\in
\bigoplus_{k=1}^n\bar H_0^1(0,\infty)\,,\
\zeta_\Phi\in\C^n\right\}\\
\equiv&\bigoplus_{k=1}^n\bar H^1(0,\infty)\,,
\end{align*}
\begin{align*}
\bar K^2:=&
\left\{\Phi=\Phi_0+\zeta_\Phi\g\,,\ \Phi_0\in\bigoplus_{k=1}^n\bar
H_0^2(0,\infty)\,,\ 
\zeta_\Phi\in\C^n\right\}\\
\equiv&\bigoplus_{k=1}^n\bar H^2(0,\infty)\,,
\end{align*}
and
$$
K^1:=\bar K^1\cap \bigoplus_{k=1}^n L^2(0,\infty)\equiv 
\bigoplus_{k=1}^nH^1(0,\infty)\,.
$$
$$
K^2:=\bar K^2\cap \bigoplus_{k=1}^n L^2(0,\infty)\equiv 
\bigoplus_{k=1}^nH^2(0,\infty)\,.
$$
One makes
$\bigoplus_{k=1}^n \left(\bar H^1(0,\infty)\oplus
L^2(0,\infty)\right)$ a Hilbert space by the scalar
product
\begin{align*}
&\langle\langle\,(\Phi,\Psi),(\tilde\Phi,\tilde\Psi)\,\rangle\rangle\\
:=&
\sum_{1\le k\le n}\langle\phi'_k,\tilde\phi'_k\rangle
+\sum_{1\le k\le n}\langle\psi_k,\tilde\psi_k\rangle+\sum_{1\le k,j\le n}\Theta_{kj}\,\bar\phi_k(0_+)\tilde\phi_j(0_+)\,.
\end{align*}
Here we put $\Phi\equiv(\phi_1,\dots,\phi_n)$,
$\Psi\equiv(\psi_1,\dots,\psi_n)$ and we used the fact that
$\zeta_\Phi=(\phi_1(0_+),\dots,\phi_n(0_+))$. \par
By Theorem 3.6 we define now skew-adjoint operators $W_\Theta$ 
corresponding to wave equations on
star-like graphs:  
the operator
$$
W_\Theta: D(W_\Theta)\to\bigoplus_{k=1}^n\left(\bar H^1(0,\infty)\oplus L^2(0,\infty)\right)\,,
$$
\begin{align*}
&D(W_\Theta):=\\&\left\{\Phi\in\bigoplus_{k=1}^n \bar H^2(0,\infty)\, 
:\, \phi_k'(0_+)+\sum_{1\le j\le n}\Theta_{k,j}\,\phi_j(0_+)
=0\,,\ 1\le k\le n\right\}\\&\oplus\, \bigoplus_{k=1}^n H^1(0,\infty)\,,
\end{align*}
$$
W_\Theta(\Phi,\Psi):=(\Psi,\,\bar A\Phi_0)\equiv \left(\Psi,\, 
\bigoplus_{k=1}^n \bar A_0\Phi_0\right)
\equiv \left(\psi_1,\dots,\psi_n,\phi''_1,\dots,\phi''_n\right)
$$
is skew-adjoint and coincides with 
$$
W: \bigoplus_{k=1}^n\left(\bar H_0^2(0,\infty)\oplus
H^1_0(0,\infty)\right)\to\bigoplus_{k=1}^n\left(\bar H_0^1(0,\infty)\oplus L^2(0,\infty)\right)\,,
$$
$$
W(\Phi,\Psi):=(\Psi,\,\bar A\Phi)\equiv \left(\Psi,\,
\bigoplus_{k=1}^n \bar A_0 \Phi\right)
\equiv \left(\psi_1,\dots,\psi_n,\phi''_1,\dots,\phi''_n\right)
$$
on the set 
$$
\left\{\Phi\in\,\bigoplus_{k=1}^n \bar H_0^2(0,\infty) :\,
\phi_k'(0_+)=0\,,\ 1\le k\le n\right\}\oplus \,\bigoplus_{k=1}^n H_0^1(0,\infty)\,.
$$ 
Moreover, by Theorem 3.7, the linear operator
\begin{align*}
&D(A_\Theta)\\:=&\left\{\Phi\in\bigoplus_{k=1}^n H^2(0,\infty)\,:\,
\phi_k'(0_+)+\sum_{1\le j\le n}\Theta_{k,j}\,\phi_j(0_+)
=0\,,\ 1\le k\le n\right\}\,,
\end{align*}
$$
A_\Theta:D(A_\Theta)\subset 
\bigoplus_{k=1}^n L^2(0,\infty)\to \bigoplus_{k=1}^n L^2(0,\infty)\,,\quad
A_\Theta\Phi:=(\phi''_1,\dots,\phi''_n)\,,
$$
is negative, injective self-adjoint, its resolvent has an
integral kernel
given by
\begin{align*}
&(-A_\Theta+\lambda^2)^{-1}(x_1,\dots,x_n,y_1,\dots,y_n)\\=&
\g_D(\lambda;x_1,y_1)\cdots\g_D(\lambda;x_1,y_1)+
\sum_{1\le k,j\le n}(\Theta+|\lambda|)^{-1}_{kj}\,e^{-|\lambda|(x_k+y_j)}\,.
\end{align*}
The operator $A_\Theta$ is of the class of Laplacian operators on
star-like graphs (see e.g. \cite{[KS]} and references therein) and the
positive quadratic form corresponding to $-A_\Theta$ is
$$
\Q:\bigoplus_{k=1}^n H^1(0,\infty)\subset 
\bigoplus_{k=1}^n L^2(0,\infty)\to\RE\,,
$$
$$
\Q(\Phi):=\sum_{1\le k\le n}\|\phi'_k\|^2_2+\sum_{1\le k,j\le
n}\Theta_{kj}\,
\bar\phi_k(0_+)\phi_j(0_+)\,.
$$
\vskip10pt\p
Let $B$ be the injective selfadjoint operator on $L^2$, the Hilbert space of square
integrable functions on $\RE^3$, 
given by $B=\sqrt{-\Delta}$. 
Then $\HU$ coincides with the Sobolev space $H^1$ of $L^2$ function
with $L^2$ distributional derivatives. $\HUB$ is nothing else that the
usual Riesz potential space $\bar H^1$ given by the set of tempered
distributions with a Fourier transform (denoted by $F$) which is
square integrable
w.r.t. the measure with density $|k|^2$. The operator $\B$ is then
defined by 
$$
F{\B\phi}\,(k):=|k|\,F\phi(k)\,.
$$   
The space $\HDB$ coincides with the space $\bar H^2$ of distributions in
$\bar H^1$ with a Fourier transform which is
square integrable w.r.t. the measure with density $|k|^2(|k|^2+1)$. 
By Sobolev embedding theorems the elements of both $\bar H^1$ and
$\bar H^2$ are ordinary functions. Indeed 
$$
\bar H^2\subset\bar H^1\subset L^6(\RE^3)\,,
\qquad \bar H^2\subset C_b(\RE^3)\,,
$$
the embeddings being continuous. 
The linear operator $\A:=-B\bar B$ acts on $\bar H^2$ as the distributional Laplacean 
$\Delta$, or equivalently
$$
F{A\phi}\,(k):=-|k|^2\,F\phi(k)\,. 
$$
In the sequel $\langle\cdot,\cdot\rangle$ will denote 
the scalar product on
$L^2$. More generally, for any $\phi$, $\varphi$ such that
$\phi\varphi$ is integrable, we will use the notation
$$\langle\phi,\varphi\rangle:=\int_{\RE^3} dx\,\bar\phi(x)\varphi(x)\,.$$
Moreover $*$ will denote convolution.
\vskip10pt\p
{\it -- Example 2.} On the Hilbert space $\bar H^1\oplus L^2$ with scalar product
$$
\langle\langle\,(\phi_1,\psi_1),(\phi_2,\psi_2)\,
\rangle\rangle:=
\langle\nabla\phi_1,\nabla\phi_2\rangle
+\langle\psi_1,\psi_2\rangle
$$
we consider the skew-adjoint operator
$$
W:\bar H^2\oplus H^1\subset \bar H^1\oplus L^2\to 
\bar H^1\oplus L^2\,,\quad W(\phi,\psi):=
(\psi,\Delta\phi)\,.
$$
by Theorem 2.5 its resolvent is given by
$$
(W+\lambda)^{-1}(\phi,\psi)=
\left(\g_\lambda*(\psi+\lambda\phi),
-\phi+\lambda \g_\lambda*(\psi+\lambda\phi)\right)\,,
$$ 
where
$$
\g_\lambda(x)=\frac{e^{-|\lambda x|}}{4\pi|x|}\,,\qquad\g\equiv\g_0\,. 
$$
Given an injective and positive Hermitean $n\times
n$ matrix $\Theta=(\theta_{ij})$, we consider the Hilbert
space $\bar H^1\oplus L^2\oplus\C^n$ with scalar product
$$
\langle\langle\,(\phi_1,\psi_1,\zeta_1),(\phi_2,\psi_2,\zeta_2)\,
\rangle\rangle_\Theta:=
\langle\nabla\phi_1,\nabla\phi_2\rangle
+\langle\psi_1,\psi_2\rangle+(\Theta
\zeta_1,\zeta_2)\,
$$
where $(\cdot,\cdot)$ denotes the scalar product on $\C^n$.\par
Given $Y=\left\{y_1,\dots,y_n\right\}\subset\RE^3$, let 
$$
\tau:\bar H^2\to\C^n\,,\quad(\tau\phi_0)_i:=\phi_0(y_i)\,,\ 1\le i\le n\,.
$$
Such a map satisfies H3.1 since $\tau^*\zeta=\zeta_i\delta_{y_i}$,
where $\delta_y$ denotes Dirac's mass at $y$. Here and below we use
Einstein's summation convention: repeated indices mean summation.\par
We define then the continuous linear map, which obviously satisfies H3.0,
$$
\tau_\Theta:\bar H^2\oplus L^2\oplus\C^n\to\C^n\,,\qquad
\tau_\Theta(\phi_0,\psi,\zeta)_i:=\phi_0(y_i)-\theta_{ij} \zeta_j\,.
$$
Thus, according to the definitions (3.2) and (3.3) one obtains
$$
\breve G_\Theta(\lambda)(\phi,\psi,\zeta)_i=
\langle\g_\lambda^{i},\psi+\lambda\phi\rangle
-\frac{1}{\lambda}
\theta_{ij} \zeta_j\,,
$$
where $\g_\lambda^{i}(x):=\g_\lambda(x-{y_i})$,
and
$$
G_\Theta(\lambda)\zeta
=\left(\lambda\zeta_i\,\g*\g^{i}_\lambda,
\,-\zeta_i\,\g^{i}_\lambda,
\,-\frac{1}{\lambda}\,\zeta\right)\,.
$$
Therefore, putting $\g^i:=\g^i_0$, by (3.4), 
\begin{align*}
&\left(\Gamma_\Theta(\lambda)\zeta\right)_i
=-\left(\tau_\Theta G_\Theta(\lambda)\zeta\right)_i\\
=&-\left(\lambda
\langle\g^{i},\g^{j}_\lambda\rangle
+\frac{1}{\lambda}\,\theta_{ij}\right)\zeta_j\\
=&-\frac{1}{\lambda}\left(\frac{|\lambda|}{4\pi}\,\zeta_i+
\left((1-\delta_{ij})\,
\frac{1-e^{-|\lambda(y_i-y_j)|}}{4\pi|y_i-y_j|}+\theta_{ij}\right)\zeta_j
\right)\,,
\end{align*}
i.e., defining
$$
\Theta_Y:=\left((1-\delta_{ij})\,
\frac{1}{4\pi|y_i-y_j|}\right)\,,\quad 
M(\lambda):=\left((1-\delta_{ij})\,
\frac{e^{-|\lambda(y_i-y_j)|}}{4\pi|y_i-y_j|}\right)\,,
$$
$$
\Gamma_\Theta(\lambda)=-\frac{1}{\lambda}
\,\left(\Theta+\Theta_Y+\frac{|\lambda|}{4\pi}-M(\lambda)\right)\ .
$$
Since H3.2 is verified by taking $\V=L^2_{\text {\rm loc}}$, we put,
defining $\g^i(x):=\g(x-y_i)$,   
$$
\bar K^1:=
\left\{\phi\in L^2_{\text {\rm loc}}\ :\ 
\phi=\phi_0+\zeta_\phi^i\g^i\,,\ \phi_0\in\bar H^1\,,\
\zeta_\phi\in\C^n\right\}\,,
$$
$$
\bar K^2:=
\left\{\phi\in L^2_{\text {\rm loc}}\ :\ 
\phi=\phi_0+\zeta_\phi^i\g^i\,,\ \phi_0\in\bar H^2\,,\
\zeta_\phi\in\C^n\right\}\,,
$$
$$
K^1:=\bar K^1\cap L^2\,,
$$
and making $\bar K^1\oplus L^2$ a Hilbert space by the scalar
product
\begin{align*}
&\langle\langle\,(\phi,\psi),(\tilde\phi,\tilde\psi)\,\rangle\rangle
_{\bar K^1\oplus L^2}\\
:=&
\langle\nabla\phi_0,\nabla\tilde\phi_0\rangle
+\langle\psi,\tilde\psi\rangle+(\Theta
\zeta_{\phi},\zeta_{\tilde\phi})\,,
\end{align*}
by Theorem 3.6 the operator
$$
D(W_\Theta):=\left\{\phi\in\bar K^2\ :\ 
\theta_{ij} \zeta_\phi^j=\phi_0(y_i)\right\}\oplus K^1\,,
$$ 
$$
W_\Theta: D(W_\Theta)\subset \bar K^1\oplus L^2\to\bar K^1\oplus
L^2\,,
$$
$$
W_\Theta(\phi,\psi):=(\psi,\,\Delta\phi_0)\equiv
\left(\psi,\,\Delta\phi+\zeta_\phi^j\,\delta_{y_j}\right)
$$
is skew-adjoint and coincides with $ W$ on the
set 
$$
\left\{\phi\in\bar H^2\, :\,  \phi(y)=0\,,\,y\in Y\right\}\oplus H^1\,.
$$ 
In the case $Y=\left\{0\right\}$ this operator coincides with the one
constructed in \cite{[BNP]}. By Theorem 3.7, the
positive quadratic form 
$$
\Q:K_1\to\RE\,,\quad\Q(\phi):=\|\nabla\phi_0\|^2_{L^2}
+\|\Theta^{1/2}\zeta_\phi\|^2_{\C^n}
$$ 
is closed and the corresponding self-adjoint operator $-\Delta_\Theta$
is defined by 
$$
D(\Delta_\Theta)=\left\{\phi\in\bar K^2\cap L^2\,:\, 
\theta_{ij} \zeta_\phi^j=\phi_0(y_i)\right\}\,,
$$
$$
\Delta_\Theta\phi:=\Delta\phi_0\,.
$$
It coincides with $\Delta$ on the set $\left\{\phi\in
H^2\,:\,\phi(y)=0\,,\ y\in Y\right\}$. Its resolvent is given by
$$
(-\Delta_\Theta+\lambda^2)^{-1}\psi=\g_\lambda*\psi+
\left(\Theta+\Theta_Y+\frac{|\lambda|}{4\pi}-M(\lambda)\right)^{-1}_{ij}
\langle\g_\lambda^i,\psi\rangle\,\g_\lambda^j\,.
$$
This operator is of the class of point perturbation of the Laplacian
(see \cite{[AGHH]} and references therein).
\vskip 10pt\p
{\it -- Example 3.} Given $v\in\RE^3$, $|v|<1$, we consider the skew-adjoint operator
$$
W^v:\bar H^2\times H^1\subseteq H^v\to H^v\,,$$
$$ W^v(\phi,\psi,z):=
(L_v\phi+\psi,L_v\psi+\Delta\phi)\,,
$$
where $L_v:=v\cdot\nabla$ and $H^v=\bar H^1\times L^2$ with scalar
product
$$
\langle\langle\,(\phi_1,\psi_1),(\phi_2,\psi_2)\,\rangle\rangle_v:=
\langle\nabla\phi_1,\nabla\phi_2\rangle+\langle L_v\phi_1,\psi_2\rangle
+\langle\psi_1,L_v\phi_2\rangle+\langle\psi_1,\psi_2\rangle\,.
$$ 
Hypotheses H4.1.1 and H4.1.2 are satisfied with $C_1=C_2=L_v$  and,
by Theorem 4.2, with $C=2L_v$ and $B=(-\Delta+L_v^2)$, 
the resolvent of $W^v$ is given by
\begin{align*}
&(W^v+\lambda)^{-1}(\phi,\psi)=
\left(\g^v_\lambda*(\psi+(-L_v+\lambda)\phi),\right.\\
&\left.-\phi+(-L_v+\lambda)\,\g^v_\lambda*(\psi+(-L_v+\lambda)\phi)\right)\,,
\end{align*} 
where
$$
F\g^v_\lambda(k)=\frac{1}{(2\pi)^{3/2}}\,\frac{1}{|k|^2+(iv\cdot k+\lambda)^2}\,. 
$$
Let
$$
\tau:\bar H^2\to\C\,,\quad\tau\phi_0:=\phi_0(0)\,.
$$
By Example 2 we know that such a map satisfies H3.1. For any real
$\theta>0$, define now the linear map, which obviously satisfies H3.0,
$$
\taub_\theta:D(\taub)\times L^2\times\C\subseteq H^v\oplus\C\to\C\,,\quad
\taub_\theta(\phi,\psi,\zeta):=\taub\phi-\theta \zeta\,,
$$
where, denoting by $\langle\phi\rangle_R$ the average of $\phi$ over
the sphere of radius $R$,
$$
D(\taub):=\left\{\phi\in\bar H^1\,:\, \lim_{R\downarrow 0}\,
\langle\phi\rangle_R\quad \text{\rm exists and is finite}\right\}\,,
\quad
\taub\phi:=\lim_{R\downarrow 0}\,
\langle\phi\rangle_R\,.
$$
Since, by Fourier transform,
$$
\taub\phi=\lim_{R\downarrow 0}\,
\frac{1}{(2\pi)^{3/2}}\int_{\RE^3}dk\,\frac{\sin R|k|}{R|k|}\,F\phi(k)\,,
$$
with reference to the notations of Section 4, we are taking here the
regularizing family 
$$J_\nu=\frac{1}{(2\pi)^{3/2}}\,\left(\nu\sqrt{-\Delta}\,\right)^{-1}{\sin
\nu\sqrt{-\Delta}}\,.$$
Thus $\bar H^2\subset D(\taub)$ and $\taub\phi=\phi(0)$ for any
$\phi\in\bar H^2$.
Then one obtains, by (4.2) and (4.4), 
$$
\breve G^v_\theta(\lambda)(\phi,\psi,\zeta)=
\langle\g^v_\lambda,\psi+(-L_v+\lambda)\phi\rangle
-\theta\lambda^{-1}\zeta
$$
and
$$
G_\theta^v(\lambda)\zeta=\zeta((-L_v+\lambda)\,\g^v*\g^v_\lambda,
-\,\g^v_\lambda,-\lambda^{-1})\,,
$$
where $$\g^v(x):=\g^v_0(x)\equiv \frac{1}{4\pi\sqrt{|x|^2-|v\wedge
x|^2}}\,.$$ 
Regarding hypothesis H4.3.1 one has
\begin{align*}
&\taub L_v\g^v*\g^v_\lambda\\
=&\frac{1}{(2\pi)^3}\,\lim_{R\downarrow 0}
\int_{\RE^3}dk\,\frac{\sin R|k|}{R|k|}\,
\frac{1}{|k|^2-(v\cdot k)^2}\,\frac{-iv\cdot
k}{|k|^2+(iv\cdot k+\lambda)^2}\\
=&\frac{1}{(2\pi)^2}\,\lim_{R\downarrow 0}
\int_0^\infty dr\,\frac{\sin Rr}{Rr}\int_0^\pi \frac{d\theta\,\sin\theta}{1-(|v|\cos\theta)^2}\,
\frac{-i|v|r\cos\theta}{r^2+(i|v|r\cos\theta+\lambda)^2}\\
=&\frac{1}{(2\pi)^2|v|}\,\lim_{R\downarrow 0}
\int_0^\infty dr\,\frac{\sin Rr}{Rr}\int_{-|v|}^{|v|} \frac{ds}{1-s^2}\,
\frac{irs}{r^2+(-irs+\lambda)^2}\\
=&\frac{1}{(2\pi)^2|v|}\,\lim_{R\downarrow 0}
\int_0^\infty dr\,\frac{\sin Rr}{Rr}\int_{-|v|}^{|v|} \frac{ds}{1-s^2}\,
\frac{-2\lambda r^2s^2}{((1-s^2)r^2+\lambda^2)^2
+4\lambda^2r^2s^2}\\
=&-\frac{4\lambda}{(2\pi)^2|v|}\,
\int_0^\infty dr\int_{0}^{|v|} ds\,\frac{s^2}{1-s^2}\,
\frac{r^2}{((1-s^2)r^2+\lambda^2)^2+4\lambda^2r^2s^2}
\,,
\end{align*}
and
\begin{align*}
&\tau \g^v*\g^v_\lambda\\
=&\frac{1}{(2\pi)^3}\,
\int_{\RE^3}dk\,
\frac{1}{|k|^2-(v\cdot k)^2}\,\frac{1}{|k|^2+(iv\cdot k+\lambda)^2}\\
=&\frac{1}{(2\pi)^2}\,
\int_0^\infty dr\,\int_0^\pi \frac{d\theta\,\sin\theta}{1-(|v|\cos\theta)^2}\,
\frac{1}{r^2+(i|v|r\cos\theta+\lambda)^2}\\
=&\frac{1}{(2\pi)^2|v|}\,
\int_0^\infty dr\,\int_{-|v|}^{|v|} \frac{ds}{1-s^2}\,
\frac{1}{r^2+(-irs+\lambda)^2}\\
=&\frac{2}{(2\pi)^2|v|}\,
\int_0^\infty dr\,\int_{0}^{|v|} \frac{ds}{1-s^2}\,
\frac{(1-s^2)r^2+\lambda^2}{((1-s^2)r^2+\lambda^2)^2
+4\lambda^2r^2s^2}
\,.
\end{align*}
Thus H4.3.1 holds true and 
\begin{align*}
&\Gamma^v_\theta(\lambda)
:=-\taub_\theta  G_\theta^v(\lambda)\\
=&-\left(\frac{\theta}{\lambda}+\frac{2\lambda}{(2\pi)^2|v|}\,
\int_0^\infty dr\,\int_{0}^{|v|} \frac{ds}{1-s^2}\,
\frac{(1+s^2)r^2+\lambda^2}{((1-s^2)r^2+\lambda^2)^2
+4\lambda^2r^2s^2}\right)\\
=&-\left(\frac{\theta}{\lambda}+\frac{\lambda}{(2\pi)^2|v|}\,
\int_{0}^{|v|} \frac{ds}{1-s^2}\,\int_{-\infty}^\infty dr\,
\frac{(1+s^2)r^2+\lambda^2}{((1-s^2)r^2+\lambda^2)^2
+4\lambda^2r^2s^2}\right)\\
=&-\frac{1}{\lambda}\,
\left({\theta}+\frac{|\lambda|}{8\pi|v|}\,
\int_{0}^{|v|}ds\, \frac{2-s^2}{(1-s^2)^2}\right)\\
=&-\frac{1}{\lambda}\,
\left({\theta}+\frac{|\lambda|}{16\pi}\,
\left(\frac{1}{1-|v|^2}+\frac{3}{2}\,\frac{1}{|v|}\,\ln{\frac{1+|v|}{1-|v|}}\,
\right)\right)
\,.
\end{align*}
Note that
$$
\lim_{|v|\downarrow 0}\,\Gamma^v_\theta(\lambda)=-\frac{1}{\lambda}\,
\left(\theta+\frac{|\lambda|}{4\pi}\right)\,,
$$
in accordance with the previous example when $Y=\left\{0\right\}$.\par
Since H3.2 is verified by 
taking $\V=L^2_{\text {\rm loc}}$, defining the Hilbert space
\begin{align*}
&H^v_\theta:=\{(\phi,\psi)\in L^2_{\text {\rm loc} }\times H^{-1}\,:\, 
\phi=\phi_0+\zeta_\phi\g^v,\\ &\psi=\psi_0-\zeta_\phi L_v\g^v, \ 
(\phi_0,\psi_0,\zeta_\phi)\in\bar H^v\oplus\C\}
\end{align*}
with scalar product
\begin{align*}
&\langle\langle\,(\phi,\psi),(\tilde\phi,\tilde\psi)\,\rangle\rangle_{H^v_\theta}
:=
\langle\nabla\phi_0,\nabla\tilde\phi_0\rangle+
\langle L_v\tilde\phi_0,\tilde\psi_0\rangle\\&
+\langle\psi_0,L_v\tilde\phi_0\rangle+\langle\psi_0,\tilde\psi_0\rangle
+\theta\zeta_{\phi}^*\zeta_{\tilde\phi}\,,
\end{align*}
by Theorem 4.11 the operator 
$$
W_\theta^v:D(W^v_\theta)\subseteq H^v_\theta\to H^v_\theta\,,
$$
\begin{align*}
D(W_\theta^v):=\{&(\phi,\psi)\in H^v_\theta\, :\,
 \phi_0=\phi_\lambda+2\zeta_\psi L_v\g^v*\g^v_\lambda
\,,\\
&\psi_0=\psi_\lambda+\zeta_\psi(\g^v_\lambda-2L_v^2\g^v*\g^v_\lambda)\,,\\
&\phi_\lambda\in\bar H^2\,,\ \psi_\lambda\in H^1\,,\  
\zeta_\psi\in\C\,,\ \theta \zeta_\phi=\taub\phi_0\,\}\,.
\end{align*}
\begin{align*}
W_\theta^v(\phi,\psi)&:=
(L_v\phi_0+\psi_0,L_v\psi_0+\Delta\phi_0)\\&\equiv
(L_v\phi+\psi,L_v\psi+\Delta\phi+\zeta_\phi\delta_0)\,,
\end{align*}
is skew-adjoint. It coincides with $W^v$ on the dense set
$$\left\{\phi\in\bar H^2\ :\
\phi(0)=0\right\}\times H^1\,. $$

\vskip 10pt\p
{\it -- A digression on the classical electrodynamics of a point particle.}
Let us begin with a discussion at the euristic level ignoring
the singular behaviour due to the self-energy of the
point particle. \par 
In the Coulomb gauge the Maxwell-Lorentz system, i.e the
nonlinear infinite dimensional dynamical system describing a
(relativistic) charged
point particle interacting with the self-generated radiation field, is
given by the equations 
\begin{align*}
\dot \As&=\Es\\
\dot \Es&=\Delta \As+4\pi e\,M\,v\delta_q\\
\dot q&=v\\
\dot p&=e\,\nabla \As(q)\,\cdot v\,,
\end{align*}
where 
$$
v=v(\As,q,p):=\frac{p-e\,\As(q)}
{\sqrt{\left|p-e\,\As(q)\right|^2+m^2}}
$$
or, equivalently,
$$
p=p(\As,q,v)=\frac{mv}{\sqrt{1-{|v|^2}}}+e\,\As(q)\,.
$$
Here we put $c=1$, where $c$ denotes the velocity of light, $e$
denotes the electric charge, $M$ is the projection onto the
divergenceless vector fields, $\As\equiv (A_1,A_2,A_3)$, div$\As=0$, 
is the vector
potential of the electromagnetic field, $q$, $v$, $|v|<1$, and $p$ denote
the particle
position, velocity and momentum respectively.
Since the total (particle + field) momentum 
$$
\Pis:=p-\frac{1}{4\pi}\,\langle \Es,\nabla \As\rangle\,,\quad 
\langle \Es,\nabla
\As\rangle:=\sum_{j=1}^{3}\int_{\RE^3}dx\,E_j(x)\nabla A_j(x)\,,
$$
is conserved, the above dynamical system can be reduced. Indeed, by defining
the fields
$$
\Phs(x):=\As(x+q)\,,\quad \Pss(x):=\Es(x+q)\,, 
$$
the Maxwell-Lorentz system can be re-written as
\begin{align*}
\dot\Phs&=v\cdot\nabla\Phs+\Pss\\
\dot\Pss&=v\cdot\nabla\Pss+\Delta\Phs+{4\pi e}\,M\,
v\delta_0\\
\dot q&=v\\
\dot\Pis&=0\,,
\end{align*}
where now 
$$
v=v(\Phs,\Pss):=\frac{p-e\,\Phs(0)}
{\sqrt{\left|p-e\,\Phs(0)\right|^2+m^2}}\,,
$$
equivalently
$$
\frac{mv}
{\sqrt{1-|v|^2}}=-e\,\Phs(0) +p\,,
$$
and
$$
p=p(\Phs,\Pss):=\Pis+\frac{1}{4\pi}\,\langle\Pss,\nabla\Phs\rangle\,.
$$
Thus we have that, at any fixed total momentum $\Pis$, we can solve the
equations for the fields $\Phs$ and $\Pss$ alone, and then recover the
particle dynamics by $\dot q= v(\Phs,\Pss)$.
\par
Due to the
singularity produced by the Dirac mass $\delta_q$, the above reasoning
is definitively  not rigorous since $\As$ is singular at the particle 
position $q$ (equivalently $\Phs$
is singular at the origin). However Example 3 suggests the definition of a 
well-defined nonlinear operator candidate to describe, in a
rigorous way, the classical electrodynamics of a point particle. \par
Let us define the infinite dimensional manifold 
\begin{align*}
{\mathcal M}:=\{(\Phs,\Pss)\,:\ &\Phs=\Phs_0+e\,M\,v\g^v,\ 
\Psi=\Pss_0-e\,M\,v\, v\cdot\nabla\g^v\,, \\ 
&(\Phs_0,\Pss_0,v)\in\bar H^1_*\times L_*^2\times \RE^3\,,\  |v|<1\}\,,
\end{align*}
where the subscript ${\,}_*$ means
``divergenceless'', $H^1$ and $L^2$ are defined as in Example 3 but
now refer to $\RE^3$-valued vector fields, and 
$$
\g^v(x):=\frac{1}{\sqrt{|x|^2-|v\wedge x|^2}}\,.
$$
Note that $$\As^v_{LW}(t,x):=e Mv\g^v(x-vt)$$
satisfies $$\square \As^v_{LW}=4\pi e Mv\delta_q\,,\quad q(t)=vt\,,$$ 
i.e. $\As^v_{LW}$ is the Li\'enard-Wiechert potential corresponding to a 
particle with constant velocity $v$.\par
We identify $T_{(\Phs,\Pss)}{\mathcal M}$, the tangent space of 
$\mathcal M$ at
$$(\Phs,\Pss)\equiv(\Phs_0+e\,M\,v\g^v,\Pss_0-e\,M\,v\,
v\cdot\nabla\g^v)\,,$$ with the Hilbert space
\begin{align*}
\H_v:=\{(\tilde\Phs,\tilde\Pss)\,:\ &\tilde \Phs=\tilde\Phs_0+e\,M\,\tilde v\g^v,\ 
\tilde\Pss=\tilde\Pss_0-e\,M\,\tilde v\, v\cdot\nabla\g^v\,, \\ 
&\tilde\Phs\in L^2_*\,,\quad(\tilde\Phs_0,\tilde\Pss_0,\tilde v)\in
\bar H^1_*\times L_*^2\times \RE^3\}\\
\simeq
\{(\tilde\Phs,\tilde\Pss)\,:\ &\tilde \Phs=\tilde\Phs_\lambda+e\,M\,
\tilde v\g^v_\lambda,\ 
\tilde\Pss=\tilde\Pss_\lambda-e\,M\,\tilde v\, v\cdot\nabla\g^v_\lambda\,, \\ 
&\quad(\tilde\Phs_\lambda,\tilde\Pss_\lambda,\tilde v)\in
H^1_*\times L_*^2\times \RE^3\}
\,.
\end{align*}
 Then we define the nonlinear vector
field 
$$X_e:D(W_e)\subset{\mathcal M}\to T{\mathcal M}\,,\quad X_e(\Phs,\Pss):=((\Phs,\Pss),W_e(\Phs,\Pss))$$ by
\begin{align*}
D(W_e):=\{&(\Phs,\Pss)\in {\mathcal M}\, :\,
 \Phs_0=\Phs_\lambda+2e\,M\,w\,v\cdot\nabla\g^v*\g^v_\lambda
\,,\\
&\Pss_0=\Pss_\lambda+e\,M\,w(\g^v_\lambda-2(v\cdot\nabla)^2\g^v*\g^v_\lambda)\,,\\
&\Phi_\lambda\in\bar H^2_*\,,\ \Psi_\lambda\in H^1_*\,,\  
w\in\RE^3\,,\\ 
&v=v(\Phi,\Pss)=\frac{p
-e\,\langle\Phs_0\rangle}
{\sqrt{\left|p
-e\,\langle\Phs_0\rangle\right|^2+m^2}}\,,\\
&p=p(\Phs,\Pss):=\Pis+\frac{1}{4\pi}\,\langle\Pss_0,\nabla\Phs_0\rangle
\,\}\,,
\end{align*}
\begin{align*}
&W_e(\Phs,\Pss):=
(v\cdot\nabla\Phs_0+\Pss_0,v\cdot\nabla\Pss_0+\Delta\Phs_0)\\
\equiv &
(v\cdot\nabla\Phs+\Pss,v\cdot\nabla\Pss+\Delta\Phs+4\pi e\,Mv\delta_0)\,.
\end{align*}
Here 
$$
F\g^v_\lambda(k):=\frac{4\pi}{(2\pi)^{3/2}}\,\frac{1}{|k|^2+(iv\cdot k+\lambda)^2}\,,
$$
$$
\langle\Phs\rangle:=\lim_{R\downarrow 0}\,\langle\Phs\rangle_R\,,
$$
with $\langle\Phs\rangle_R$ denoting the average of $\Phs$ over the
sphere of radius $R$, and $\lambda$ is an arbitrary positive
parameter. We remark that, as it should be clear from the general results
given in the previous sections, the parameter $\lambda$ has simply the
role of allowing a convenient decomposition (into ``regular'' and
``singular'' components) of the elements in $D(W_e)$, but plays no
role in the definition of the action of $W_e$, which
indeed is $\lambda$-independent.
\par
It is not difficult to check, by a direct computation, that 
$$W_e(\Phs,\Pss)\in \H_{v(\Phs,\Pss)}\,,$$ 
so that $X_e$ is a vector field on ${\mathcal M}$ in the
differential geometric sense as stated above.\par
Note that $W_e$ coincides with the linear operator $W_0$ corresponding
to the free wave equation on 
the dense set $\left\{\Phs\in\bar H_*^2\
:\ \Phs(0)=0\right\}\times H^1_*$, so that $W_e$ is a nonlinear
singular perturbation of the skew-adjoint $W_0$. 
\par
Once the vector field $X_e$ is defined, the first question to be posed is: 
does $X_e$ generate a nonlinear flow $F_e(t)\,$? At the present we have no definitive answer to this question. The
results obtained in the linear case (see \cite{[NP]}) suggest to try
to write the presumed solution as 
$F_e(t)((\Phs(0),\Pss(0)))\equiv (\Phs(t),\Pss(t))=(\Phs_v(t),\Pss_v(t))$
where $v=v(t)$ is a pre-assigned time-dependent vector and 
$(\Phs_v(t),\Pss_v(t))$ is the solution of the linear inhomogeneous,
time-dependent, wave equation
\begin{align*}
\dot\Phs_v(t)&=v(t)\cdot\nabla\Phs_v(t)+\Pss_v(t)\\
\dot\Pss_v(t)&=v(t)\cdot\nabla\Pss_v(t)+\Delta\Phs_v(t)+4\pi eM\,v(t)\,\delta_0\\
\end{align*}
with initial data $(\Phs(0),\Pss(0))\in D(W_e)$, and then looking for
the right differential equation to be satisfied by $v(t)$ in order
that the fields $\Phs_v(t)$ and $\Pss_v(t)$ belong to $D(W_e)$ for any
$t$ (and hence
fit the correct nonlinear boundary conditions).


\section{Appendix: Skew-adjoint extensions of skew-symmetric operators}
Let $$W:D(W)\subseteq\H\to\H$$ be a skew-adjoint operator on the Hilbert
space $\H$ with scalar product $\langle\cdot,\cdot\rangle$ and
corresponding norm $\|\cdot\|$. The linear
subspace $D(W)$
inherits a Banach space structure by introducing the graph norm
$$\|\phi\|^2_W:=\|\phi\|^2+\|W\phi\|^2\,.$$ 
Thus, for any $\lambda\in\RE$, $\lambda\not= 0$, 
$(-W+\lambda)^{-1}\in \BO(\H,D(W))$.\par 
We consider now a linear operator
$$
L:D(L)\subseteq\H\to\fh \,,
$$
$\fh$ a Hilbert space, such that:\p
A.1)
$$
D(W)\subseteq D(L)\quad\text{\rm and}\quad
L_0:=L_{\,\left| D(W)\right.}\in \BO(D(W),\fh)\,;
$$
A.2)
$$
\ran(L_0)=\fh\,;
$$
A.3)
$$
\overline{\ker(L_0)}=\H\,.
$$
By A.3 $W_0:=W_{\,\left| \ker(L_0)\right.}$  
is a closed densely defined skew-symmetric
operator. We want now to define a skew-adjoint extension $\hat W\not=
W$ of $W_0$. It will be a singular perturbation of $W$ since it will
differ from $W$ only on the complement of the dense set $\ker(L_0)$.\par
We define, for any $\lambda\in\RE$, $\lambda\not= 0$, the 
following bounded operators:
$$
\breve G(\lambda):=L (-W+\lambda)^{-1} :\H\to\fh\,,
\qquad
G(\lambda):=-\breve G(-\lambda)^*:\fh\to\H\, .
$$
By the first resolvent identity one easily obtains the following (see
\cite{[P1]}, Lemma 2.1)
\begin{lemma} For any $\lambda\not =0$ and $\mu\not=0$ one has
\begin{eqnarray*}
(\lambda-\mu)\,\breve G(\mu)(-W+\lambda)^{-1}&
=&\breve G(\mu)-\breve G(\lambda)\\
(\lambda-\mu)\,(-W+\mu)^{-1}G(\lambda)&
=&G(\mu)- G(\lambda)\ .
\end{eqnarray*}
\end{lemma}
By \cite{[P3]}, Lemma 2.1, and A.2 one has that A.3 is equivalent to\p
A.3)
$$
\ran(G(\lambda))\cap D(W)=\left\{0\right\}\,.
$$
We further suppose that\p
A.4)
$$
\ran(G(\lambda))\subseteq D(L)\quad{\text{\rm and}}\quad L\, G(\lambda)\in\BO(\fh)\,.
$$ 
Thus we can define $\Gamma(\lambda)\in\BO(\fh)$ by 
$$
\Gamma(\lambda):=-L G(\lambda)\,
$$
and we suppose that\p
A.5)
$$
\Gamma(\lambda)^*=-\Gamma(-\lambda)\,.
$$
By lemma 6.1 one has
\begin{equation}
\Gamma(\lambda)-\Gamma(\mu)=-L_0 (G(\lambda)-G(\mu))=
(\lambda-\mu)\,\breve G(\mu)G(\lambda)
\end{equation}
and thus, by A.5 and \cite{[P1]}, Proposition 2.1, 
the operator $\Gamma(\lambda)$ is boundedly invertible for any real 
$\lambda\not =0$. 
\begin{theorem} For any real $\lambda\not =0$, under the hypotheses
{\rm A.1-A.5}, the bounded linear operator 
$$
(-W+\lambda)^{-1}+G(\lambda)\Gamma(\lambda)^{-1}\breve G(\lambda)
$$
is a resolvent of a skew-adjoint
operator $\hat W$ such that 
$$
\left\{\phi\in D(\hat W)\cap D(W)\,:\, \hat W\phi=W\phi\right\}=\ker(L_0)\,.
$$
It is defined by
$$
D(\hat W):=\left\{\,\phi\in\H\, :\, \phi=
\phi_\lambda+G(\lambda)\Gamma(\lambda)^{-1}L_0\,\phi_\lambda,
\quad \phi_\lambda\in D(W)\,\right\}\, ,
$$
$$
(-\hat W+\lambda)\phi:=(-W+\lambda)\phi_\lambda\, .
$$
Such a definition is $\lambda$-independent and the 
decomposition of $\phi$ entering
in the definition of the domain is unique. 
\end{theorem}
\begin{proof} By (6.1) 
$\hat R(\lambda):=
(-W+\lambda)^{-1}+G(\lambda)\Gamma(\lambda)^{-1}\breve G(\lambda)$
satisfies the resolvent identity $(\lambda-\mu)\,\hat R(\mu)\hat R(\lambda)=
\hat R(\mu)-\hat R(\lambda)$ (see \cite{[P1]}, page 115, for the explicit
computation) and, by A.5, 
$\hat R(\lambda)^*=-\hat R(-\lambda)$. Moreover, by A.3, 
$\hat R(\lambda)$ is injective. Thus 
$\hat W:=-\hat R(\lambda)^{-1}+\lambda$ is well-defined on 
$D(\hat W):=\ran( R(\lambda))$, is $\lambda$-independent and is
skew-symmetric. It is skew-adjoint since $\ran(W\pm\lambda)=\H$ by 
construction.
\end{proof}

\end{document}